\documentclass[12pt]{article}

\usepackage{amsmath}
\usepackage{t-angles}

\title{
Duality Theorem  and Drinfeld Double  in Braided Tensor Categories
 \thanks {This work is supported by National Science Foundation  (No. 19971074)}}
\author{
Shouchuan Zhang \\ Department  of Mathematics, Hunan  University,
Changsha  410082, P.R.China\\
E-mail:z9491@yahoo.com.cn }
\date{}
\begin{document}
\newtheorem{Theorem}{\quad Theorem}[section]
\newtheorem{Proposition}[Theorem]{\quad Proposition}
\newtheorem{Definition}[Theorem]{\quad Definition}
\newtheorem{Corollary}[Theorem]{\quad Corollary}
\newtheorem{Lemma}[Theorem]{\quad Lemma}
\newtheorem{Example}[Theorem]{\quad Example}
\maketitle \addtocounter{section}{0}

 \begin {abstract} Let $H$ be a finite Hopf algebra with  $C_{H,H} = C_{H,H}^{-1}.$
  The duality  theorem is shown for $H$, i.e., $$   (R \# H)\# H^{\hat *}   \cong R \otimes (H \bar
  \otimes H^{\hat *} ) \hbox { \ \ \  as algebras in } {\cal C}.$$
Also, it is proved that the Drinfeld double $(D(H),[b])$ is a
quasi-triangular Hopf
algebra in ${\cal C}$.
 \end {abstract}

\section {Introduction and Preliminaries}

               It is well-known that in the work of $C^*$-algebras and
               von Neumann algebras, for an abelian group $G$,the product of $R$ crossed
by $G$ crossed by $\hat G$  is isomorphic to the tensor product of
$R$ and  the  compact operator.
 Its generalization to Hopf-von Neumann algebras was known again.( see,  for example,
  Stratila \cite {St81}). Blattner and Montgomery
strip off the functional analysis and  duplicate the result at the
level of Hopf algebras (see \cite {BM85} and \cite {Mo93}). They
proved that for an ordinary Hopf algebra $H$ and some subalgebra
$U$ of $H^*$,
\begin {eqnarray*} \label {e1}
(R \# H)\# U   \cong R \otimes (H \# U ) \hbox { \ \ \ as algebras
\ . }
 \end {eqnarray*}

The basic construction of the Drinfeld double is due to Drinfeld
\cite {Dr86}. S.Majid  \cite {Ma90a} and D.E.Radiford \cite
{Ra93a} modified the treatment.

In this paper, we generalize the duality theorem and Drinfeld
double into the braided case, i.e., for a finite  Hopf algebra $H$
with  $C_{H,H} = C_{H,H}^{-1},$ we show that
  \begin {eqnarray*} \label {e2}
(R \# H)\# H^{\hat *}   \cong R \otimes (H \bar
  \otimes H^{\hat *} )  \hbox { \ \ \ as
algebras in } {\cal C}
 \end {eqnarray*}
We also show that the  Drinfeld double (D(H),[b]) is a
quasi-triangular Hopf algebra in ${\cal C}$.

In this paper, $({\cal C}, \otimes, I,  C )$  is always   a
braided tensor category, where $I$ is the identity object and $C$
is the braiding.

By \cite [Theorem 0.1] {ZC},
 we can view every braided tensor category as
 a strict braided tensor category.

  For an object $V $ in  ${\cal C}$, if there exists
an object $U$ and morphisms : $d_V : U\otimes V \rightarrow I
\hbox { \ \ \ and \ \ \ } b_V : I \rightarrow V \otimes U$ in
${\cal C}$ such that $ (d_V \otimes id _U) (id _U \otimes b_V)
=id_ U \hbox {\ \ \ and \ \ \ } (id_V \otimes d _V)(b_V \otimes id
_V)= id_ V,$ then $U$  is called  a  left dual of $V,$ written as
$V^{ *}.$ In this case, $d_V$ and $b_V$ are called the evaluation
morphism and coevaluation morphism of $V$, respectively. In
general, we use $d$ and $b$ instead of $d_V$ and $b_V.$
Furthermore, $V$ is said to be finite if $V$ has a left dual (see
\cite {Ta} ).

Let us  define the transpose $f^* =(d \otimes id _{U^*} )(id_{V^*}
\otimes f \otimes id _{U^*}) (id _{V^*}\otimes b): V^* \rightarrow
U^*$ of
 a morphism $f : U \rightarrow V$.

Let $\Delta ^{cop} = C_{H,H}\Delta$  and $m^{op} = m C_{H,H}.$ We
denote $(H, \Delta ^{cop}, \epsilon)$  by $H^{cop}$ and  $(H,
m^{op}, \eta)$  by $H^{op}$. Furthermore we denote  $((H^{ * })^{\
op }) ^ {  cop }$ by $H^{\hat *}$.

A bialgebra $(H, R, \bar \Delta )$ with convolution-invertible $R$
in $Hom _{\cal C} (I, H\otimes H)$ is called a quasi-triangular
bialgebra
 in ${\cal C}$ if
$(H, \bar \Delta , \epsilon )$  is a coalgebra and the following
conditions are satisfied:

 (QT1)  $(\bar \Delta \otimes id ) R = (id \otimes id  \otimes m )
(id \otimes R  \otimes id)R$;

 (QT2)  $(id   \otimes \Delta ) R = (m  \otimes id  \otimes id)
(id \otimes R  \otimes id)R$;

 (QT3)  $(m \otimes m )(id \otimes C_{H,H} \otimes id  )(\bar \Delta \otimes R )
=  (m \otimes m )(id \otimes C_{H,H} \otimes id  )(R \otimes
\Delta ).$

 \noindent In this case, we also
say that $(H,R, \bar \Delta )$ is a braided quasi-triangular
bialgebra.

Dually, we can define a coquasi-triangular bialgebra $(H, r, \bar
m)$  in the  braided tensor category  ${\cal C}$ .

In particular, we  say that  $(H,R)$  is quasi-triangular if $(H,
R, \Delta ^{cop})$ is quasi-triangular. Dually, we  say that
$(H,r)$  is coquasi-triangular if $(H, r, m ^{op})$ is
coquasi-triangular.

A morphism $\tau$ from $H\otimes A$ to $I$ in ${\cal C}$ is called
a skew pairing on $H\otimes A$ if the following conditions are
satisfied:

(SP1): $\tau (m \otimes id ) = \tau (id \otimes \tau \otimes id )
(id \otimes id  \otimes \Delta )$ ;

(SP2): $\tau ( id \otimes m ) = (\tau  \otimes \tau )(id \otimes
C_{H,A} \otimes id ) (\Delta \otimes id \otimes id)$;

(SP3):$\tau (id \otimes \eta ) = \epsilon _H$;

(SP4):$\tau (\eta \otimes id ) = \epsilon _A$.

 If a morphism $\tau$ from $V\otimes W$ to $I$ in ${\cal C}$
 satisfies
 $$(id _U\otimes \tau ) (C_{V,U}\otimes id_W )= (\tau \otimes id _U  ) (id_V\otimes C_{U,W} )$$
for $U = V, W,$ then $\tau $ is said  to be symmetric with respect
to the braiding $C$ \ .

\begin {Lemma} \label {01''} If $H$ has a left dual $H^*$, then the
 following conditions are equivalent:

(i) The evaluation of $H$ is symmetric with respect to the
braiding $C$.

(ii) $C_{U,V} = C_{V,U}^{-1}$  for $U, V = H $ or $H^*$.

(iii)  $C_{H,H} = C_{H,H}^{-1}$.

(iv)  $C_{H^*, H^*} = C_{H^*, H^*}^{-1}$.

(v) $(id _H\otimes d ) (C_{H^{ *},H}\otimes id_H )= (d \otimes id
_H  ) (id_{H^{ *}}\otimes C_{H,H} )$.

(vi) $(id _{H ^{ *}} \otimes d ) (C_{H^{ *},H^{ * }}\otimes id_H
)= (d \otimes id _{H^{ * }}  ) (id_{H^{ *}}\otimes C_{H^{ *},H}
)$.

 \end {Lemma}

{\bf Proof.} It is straightforward.\ \
\begin{picture}(5,5)
\put(0,0){\line(0,1){5}}\put(5,5){\line(0,-1){5}}
\put(0,0){\line(1,0){5}}\put(5,5){\line(-1,0){5}}
\end{picture}\\
\section {   The Duality Theorem }

                            In this section, we obtain the duality
theorem
 for Hopf algebras living in the braided tensor category ${\cal C}.$

Throughout this section, $H$  is a finite Hopf algebra  with
$C_{H,H}= C_{H,H}^{-1}$ living in ${\cal C}$.
 $(R, \alpha )$  is a left $H$-module algebra in ${\cal C}$ and
 $R\# H$  is the smash product in ${\cal C}$. Let $H^{\hat * }= ((H^ {  * })^{ op })^{  cop }. $

The proof of   Lemmas \ref {1.1}
---\ref {1.7} is very similar to that of correspending results in \cite {Mo93}.
\begin {Lemma} \label {1.1}
If $(R, \alpha )$ is an $H$-module algebra, let $\phi = (id
\otimes \alpha )(b' \otimes id ):$ \ $R \rightarrow  H^{*\ op }
\otimes R,$ where $b' = C_{H, H^{* \ op }}b$. Then
$$\phi m =
(m\otimes m )(id \otimes C_{R, H^{*\ op}}\otimes id)(\phi \otimes
\phi ). $$
\end {Lemma}

\begin {Lemma} \label {1.2}
(1)  $(H, \rightharpoonup )$  is a left  $H^{ \hat *}$-module
algebra under the module operation  $\rightharpoonup \  = (id
\otimes d) (C_{H^{\hat *}, H} \otimes id  ) (id \otimes \Delta )
$.

(2)  $(H^{\hat *}, \rightharpoonup )$  is a left $H$-module
algebra under the module operation  $\rightharpoonup \ = (id
\otimes d)(id \otimes C_{H, H^{\hat *}}) (C_{ H,H^{\hat *}}
\otimes id  ) (id \otimes \Delta )$.

(3)  $(H, \leftharpoonup )$  is a right $H^{ \hat *}$-module
algebra   under the module operation  $\leftharpoonup \ = (
d\otimes id) (C_{H,H^{\hat *} } \otimes id  )(id \otimes C_{H,
H^{\hat *}}) (\Delta \otimes id )$.

(4)  $(H^{\hat *}, \leftharpoonup )$  is a right $H$-module
algebra under the module operation  $\leftharpoonup \  = (
d\otimes id) (id \otimes C_{ H^{\hat *},H}) (\Delta \otimes id )$.

\end {Lemma}

\begin {Lemma}  \label {1.3} The object $H \otimes H^{\hat *}$ becomes an algebra, written as   $H
\bar \otimes H^{\hat *}$, under the  multiplication $m_{H \bar
\otimes H^{\hat *}} = (id _ H \otimes d \otimes id _{H^{\hat *}}
)$  and unity $\eta _{H \bar \otimes H^{\hat *}} = b. $
          \end {Lemma}

          In fact, if ${\cal C}$  is a braided  tensor category
          determined  by the (co)quasi-triangular structure of
          a (co)quasi-triangular Hopf algebra over a field $k$,
           then $H\bar \otimes
          H^{\hat *}$  can be viewed as  $End _k(H)$ or $M_n (k)= \{A \mid A \hbox { is an  } n\times n
\hbox { matrix  over the field } k\}$.

\begin {Lemma} \label {1.4}
Let $\lambda = (m \otimes d \otimes id ) (id \otimes C_{H^{\hat
* }, H}\otimes id \otimes id )(id \otimes id \otimes \Delta \otimes id )(id \otimes id \otimes
b)$  and $\rho  =  ( d \otimes m\otimes id ) (id \otimes id
\otimes C_{ H, H}\otimes id  )(id \otimes C_{H, H}\otimes id
\otimes id )(id \otimes id \otimes \Delta \otimes id )(id \otimes
id \otimes b)$. Then $\lambda $ is an algebra morphism from $H\#H
^{\hat *}$ to $H\bar \otimes H ^{\hat *}$ and $\rho $ is an
anti-algebra morphism from $H ^{\hat *}\#H$ to $H\bar \otimes H
^{\hat *}$.
\end {Lemma}

\begin {Lemma} \label {1.5} The following relation holds:
 $ m (\lambda
\otimes \rho ) =  m (\rho \otimes \lambda )(id \otimes
\rightharpoonup \otimes \leftharpoonup \otimes id )(id \otimes
C_{H, H^{\hat *}} \otimes C_{H^{\hat *},H} \otimes id) (id \otimes
id \otimes C_{H^{\hat *}, H^{\hat *}} \otimes id \otimes id )(id
\otimes C_{H^{\hat *},H} \otimes C_{H, H^{\hat *}} \otimes id)(
 C_{H^{\hat *},H^{\hat *}}\otimes  id \otimes id \otimes C_{H^{\hat *}, H^{\hat *}}
 )(S \otimes id  \otimes id  \otimes id  \otimes id  \otimes id )(\Delta  \otimes id
  \otimes id  \otimes \Delta  ) (id \otimes C_{H,H}\otimes id ) (C_{H, H^{\hat *} } \otimes
  C_{ H^{\hat *},H })(id \otimes C_{H ^{\hat *},H ^{\hat *}}\otimes id
  )$.
\end {Lemma}
\begin {Lemma} \label {1.6}
If the antipode of $H$ is  invertible, then  $\lambda $ is invertible.
\end {Lemma}

\begin {Lemma} \label {1.7}
$R \# H$ becomes    an  $H^{ *}$-module algebra under the
 module operation  $\rightharpoonup ' \ =  (id \otimes
\rightharpoonup )(C_{H^{\hat *}, R} \otimes id ).$
\end {Lemma}

          \begin {Theorem} \label {1.8}
 If $H$ is a  finite Hopf algebra with  $C_{H,H}=C_{H,H}^{-1}$,
then
  $$   (R \# H)\# H^{\hat *}   \cong R \otimes (H \bar
  \otimes H^{\hat *} ) \hbox { \ \ \  as algebras in }  {\cal C},$$

Where $H \bar
  \otimes H^{\hat *} $ is defined in Lemma \ref {1.3}.

\end {Theorem}

{\bf Proof. }   We first define a morphism  $w$ from $H^{\hat *} $
to $H \# H^{\hat *}$ such that  $w = \lambda ^{-1}\rho (S^{-1}
\otimes \eta _H )$; this can be done by Lemma \ref  {1.6} and
\cite [Theorem 4.1]{Ta} . Since $\rho $ and $S^{-1}$ are
anti-algebra morphisms, $w$ is an algebra morphism.  Set $ \phi =
(id \otimes \alpha  ) (C_{H, H^{\hat *}} \otimes id )(b \otimes id
)$.

We now define a morphism  $\Phi = (id \otimes  m_{H \# H^{\hat
*}})(id  \otimes w \otimes id \otimes id )(C_{ H^{\hat *}, R
}\otimes id \otimes id ) (\phi \otimes id \otimes id )$ from
$(R\#H)\# H^{\hat *}$ to $R \otimes (H \bar \otimes H^{\hat *})$
and a morphism
 $\Psi = (id \otimes  m_{H \#
H^{\hat *}})(id  \otimes w \otimes id \otimes id )(C_{ H^{\hat *
},R}\otimes id \otimes id ) (S \otimes id \otimes id \otimes id
)(\phi \otimes id \otimes id )$ from $R \otimes (H \bar \otimes
H^{\hat *})$ to $(R\#H)\# H^{\hat *}.$
 It is straightforward to verify  that $\Phi \Psi = id  $ and $\Psi \Phi =
id$.

To see that $\Phi$ is an algebra morphism, we only need to show
that $\Phi ' = (id \otimes \lambda )\Phi $ is an algebra morphism.
Set  $\xi = (id \otimes \rho  )(S^{-1} \otimes \eta _H) C_{H^{\hat
*}, R} \phi $ , which is a morphism from $R$ to $R \otimes (H \bar
\otimes H^{\hat *} )$. We have that $\xi $ is an algebra morphism
and  $\Phi ' = (id \otimes m) (\xi \otimes \lambda ).$ Using Lemma
\ref {1.5}, we can show that $ (id \otimes m)(C_{H \bar \otimes
H^{\hat *}, R}\otimes id ) (\lambda \otimes \xi) = (id \otimes m)
( \xi \otimes\lambda ) (\alpha \otimes id \otimes id  ) (id
\otimes C_{H, R} \otimes id )(id \otimes id \otimes C_{H^{\hat *},
R} )(\Delta \otimes id \otimes id ).$ Applying this, we see that
$\Phi'$ is an algebra morphism.\ \
\begin{picture}(5,5)
\put(0,0){\line(0,1){5}}\put(5,5){\line(0,-1){5}}
\put(0,0){\line(1,0){5}}\put(5,5){\line(-1,0){5}}
\end{picture}\\

 \section {Drinfeld Double}

In this section, we  construct the
 Drinfeld double $D(H)$ for  a finite Hopf algebra $H$ with  $C_{H,H} = C_{H,H} ^{-1}$
in the braided tensor category ${\cal C}$. We show that
$(D(H),[b])$ is quasi-triangular.

\begin {Theorem} \label {2.1}
Let $H$
and $A$ be two bialgebras  in ${\cal C}$. Assume that  $\tau$  is  an invertible
skew pairing on $H \otimes A$ and  symmetric  with respect to the braiding. If we
 define
$\alpha = (\tau \otimes id \otimes \bar \tau ) ( id \otimes id
\otimes  C_{H, A} \otimes id )( id \otimes   C_{H, A} \otimes
\Delta  )  (\Delta \otimes \Delta )$ and $\beta = (\tau \otimes id
\otimes \bar \tau ) ( id \otimes  C_{H, A} \otimes id \otimes id
)(\Delta \otimes   C_{H, A} \otimes   id )  (\Delta \otimes \Delta
),$
 then  the double cross product $A _\alpha \bowtie _\beta H$, defined in \cite [P36]{ZC}),
  of $A$ and $H$  is an algebra and a
coalgebra.
 If $A$ and $H$ are Hopf algebras, then $A _\alpha \bowtie _\beta H$ has an antipode.
 Furthermore, $A _\alpha \bowtie _\beta H$ is a bialgebra if and only
 if
$ C_{A, H}C_{H, A} = id .$
\end {Theorem}
{\bf Proof.} We can check that $(A, \alpha )$  is an $H$-module
coalgebra and $(H, \beta )$ is an $A$-module coalgebra step by
step. We can also  check that (M1)--(M3) in \cite [P36--37] {ZC}
hold  step by step. Consequently,  it follows from \cite
[Corollary 1.8 , Theorem 1.5]{ZC} that  $A _\alpha \bowtie _\beta
H$ is an algebra and a coalgebra. Observing  the proof of \cite
[Theorem 1.5] {ZC}, we know that the condition (M4) is not needed
in the proof.  Consequently, $A _\alpha \bowtie _\beta H$ has an
antipode.

From \cite [Proposition 3.6] {BD}, we obtain our last assertion.\
\ \
\begin{picture}(5,5)
\put(0,0){\line(0,1){5}}\put(5,5){\line(0,-1){5}}
\put(0,0){\line(1,0){5}}\put(5,5){\line(-1,0){5}}
\end{picture}\\

In this case, $A _\alpha \bowtie _\beta H$ can be
written as $A  \bowtie _\tau H$ and  called a double cross
product.

\begin {Theorem} \label {2.2}
Let $H$ be a finite  Hopf algebra with $C_{H,H}= C_{H,H}^{-1}$.
Set $A= (H^*)^{op}$ and $ \tau = d_H C_{H,A}$. Then $(D(H), [b])$
is a quasi-triangular
 Hopf algebra in ${\cal C}$ with  $[b]= \eta _A \otimes b \otimes \eta
_H$ and $D(H) =A \bowtie _{\tau} H,$ called the Drinfeld double of
$H$.
\end {Theorem}
{\bf Proof.} Using \cite [Proposition 2.4] {Ma95a} or the
definition of the evaluation and coevaluation on tensor product,
we can obtain that
 $\tau $ is a skew pairing on $H \otimes A$ and
 $[b]$ satisfies (QT1) and  (QT2). For  (QT3), see that
$(m \otimes m )(id \otimes C_{D(H),D(H)} \otimes id  )( \Delta
^{cop}\otimes [b] ) =(id \otimes id  \otimes m \otimes id ) (id
\otimes C_{A,H} \otimes id \otimes id )(id \otimes id  \otimes m
\times id  \otimes id )(id \otimes id  \otimes id \otimes m \times
id \otimes id )(id \otimes id \otimes id \otimes id  \otimes m
\times id \otimes id ) (C_{A,A} \otimes id \otimes id \otimes
C_{H,H} \times C_{H,A} )(id \otimes id \otimes id \otimes C_{H,H}
\otimes C_{H,H}\otimes  id ) (id \otimes id
 \otimes C_{H,H} \otimes C_{H,H}\otimes b )
(id \otimes id \otimes id
 \otimes C_{H,H} \otimes S^{-1} )
 (\Delta  \otimes \Delta
 \otimes C_{H,H}   )
(id   \otimes \Delta
 \otimes id )
(id   \otimes \Delta ) = (id \otimes C_{A,H} \otimes id  ) (id
\otimes m  \otimes id \otimes id  )
 (C_{A,A} \otimes C_{H,A}  \otimes id   )
  (id \otimes id \otimes  m  \otimes id  \otimes id )
 (id \otimes id   \otimes C_{H,H} \otimes C_{H,A}  )
 (id \otimes id   \otimes id \otimes C_{H,H} \otimes id  )
 (\Delta \otimes \Delta   \otimes b )$
 and $ (m \otimes m )(id \otimes C_{D(H),D(H)} \otimes id  )([b]
\otimes \Delta )= (id \otimes m  \otimes id \otimes id ) (C_{H,A}
\otimes C_{A,H}  \otimes id  ) (id \otimes id \otimes  m  \otimes
id \otimes id ) (id \otimes id \otimes m  \otimes id \otimes id
\otimes id  )(id \otimes id \otimes m  \otimes id \otimes id
\otimes id  \otimes id)(id \otimes C_{A,A} \otimes id \otimes  S
\otimes id \otimes id \otimes id  ) (id \otimes id \otimes C_{A,A}
\otimes id \otimes id \otimes id \otimes id  ) (id \otimes C_{A,A}
\otimes \Delta \otimes id \otimes id \otimes id  )(b \otimes
\Delta \otimes id \otimes id \otimes id  )(\Delta \otimes \Delta
)= (id \otimes C_{A,H} \otimes id  ) (id \otimes m  \otimes id
\otimes id )
 (C_{A,A} \otimes C_{H,A}  \otimes id   )
  (id \otimes id \otimes  m  \otimes id \otimes id )
 (id \otimes id   \otimes C_{H,H} \otimes C_{H,A}  )
 (id \otimes id   \otimes id \otimes C_{H,H} \otimes id  )
 (\Delta \otimes \Delta   \otimes b )$. Thus
 (QT3) holds.\ \
\begin{picture}(5,5)
\put(0,0){\line(0,1){5}}\put(5,5){\line(0,-1){5}}
\put(0,0){\line(1,0){5}}\put(5,5){\line(-1,0){5}}
\end{picture}\\

 In fact, there exists a very closed relation between the
Drinfeld double $D(H)$ defined in Theorem \ref {2.2} and the
Drinfeld double  ${\cal D} (H)$ defined in \cite {Ch}. Since
$C_{H^*, H} : \ \ D(H)= H^{*op} \bowtie _{\tau} H \longrightarrow
{\cal D} (H^{op})= H^{op} \bowtie _d H^* $ is an anti-algebra
isomorphism  and a coalgebra isomorphism, we have  $D(H) \cong
{\cal D}(H^{op})^{op} \ \  \hbox{ as Hopf algebras in } {\cal C}.
$

\section { Example }

In this section, using preceding conclusions, we give some
examples for the duality theorem and Drinfeld double in a braided
tensor category ${\cal C}.$

\begin {Proposition} \label {3.1}(See \cite [Definition 2.8 (R2)] {LZ})
  A $\chi$-Hopf algebra $H$ is a Hopf algebra living in the braided tensor
category $({}^{kG}{\cal M}, C^r)$, where $r(g,h) = v ^ {\chi (g,
h)}$ for any $g, h \in G$; $\chi$ is a map from $G \times G$ to
${\bf Z}$ with $\chi (a,b) = \chi (b,a)$ and $\chi (a +b, c) =
\chi(a, c)+ \chi (b, c)$ for any $a, b , c \in G$.
\end {Proposition}

It follows from the preceding proposition and \cite [Ex. 2.12,
Section 3] {LZ} that
 Lusztig's algebra  ${\ }'f$, the twisted Ringel-Hall algebra and Ringel's
 composition algebra  are all
Hopf algebras in braided tensor categories.

\begin {Example} \label {3.2}  The bilinear map $\tau$, defined in
\cite [Pro. 1.2.3]{L}, of
 Lusztig's algebra ${\ }'f$
is  symmetric with respect to the braiding.
\end {Example}
{\bf Proof.} For any homogeneous elements  $x, y, z \in { }'f$, we
have that $(\tau \otimes id)(id \otimes C_{{\ }'f,{\ } 'f} ) (x
\otimes y \otimes z ) = \tau (x,z)v^{\mid y \mid \cdot \mid z\mid
} \delta _{\mid x\mid , \mid z \mid} = \tau (x,z)v^{\mid y \mid
\cdot \mid x\mid } \delta _{\mid x\mid , \mid z \mid} = (id\otimes
\tau )( C_{{\ }'f, {\ }'f} \otimes id ) (x \otimes y \otimes z ).$
\ \
\begin{picture}(5,5)
\put(0,0){\line(0,1){5}}\put(5,5){\line(0,-1){5}}
\put(0,0){\line(1,0){5}}\put(5,5){\line(-1,0){5}}
\end{picture}\\

\begin {Example} \label {3.3} (see \cite [Example 9.4.9]{Ma95b}) The evaluation of the braided group analogue
 $\underline H$ of an ordinary coquasi-triangular cocommutative Hopf algebra $(H,r)$ is
 symmetric  with respect to braiding $C^r$. In particular, the above conclusion holds for  $H =kG$.
\end {Example}
{\bf Proof. } It is straightforward.\ \
\begin{picture}(5,5)
\put(0,0){\line(0,1){5}}\put(5,5){\line(0,-1){5}}
\put(0,0){\line(1,0){5}}\put(5,5){\line(-1,0){5}}
\end{picture}\\

\begin {Example} \label {3.4} Let $H$ denote  Lusztig's algebra ${ }'f$.
If  $A = {}'f ^{op}$,  then the bilinear map $\tau $ as in Example
\ref {3.2} is a skew pairing on $H \otimes A$ and symmetric with
respect to the braiding. Thus, by Theorem \ref {2.1}, $A  \bowtie
_\tau H$ is an algebra and a coalgebra with an antipode in $({}^{k
{\bf Z}[I]}{\cal M}, C^r) $, but it  is never a bialgebra in
$({}^{k {\bf Z}[I]}{\cal M},C^r) $ since $ C_{A, H}C_{H, A} \neq
id .$
\end {Example}

It has been known that the category of comodules of  every
ordinary coquasi-triangular Hopf algebra is a braided tensor
category. For example, let $H={\bf C Z}_n$ and $r(a,b) = (e
^{\frac {2 \pi i}{n}})^{ab}$  for any $a, b \in {\bf Z}_n$, where
${\bf C}$ is the complex field. It is clear that $({\bf CZ}_n, r
)$  is a coquasi-triangular Hopf algebra. Thus, $ ({} ^{{\bf CZ}
_n} {\cal M} , C^r)$ is a braided tensor category, usually written
as ${\cal C}_n.$  Every algebra or Hopf algebra living in ${\cal
C}_n $ is called an anyonic algebra or anyonic Hopf algebra (see
\cite [Example 9.2.4] {Ma95b}).
 Every algebra or Hopf algebra living   in ${\cal C}_2 $
is called a superalgebra or super-Hopf algebra. In particular,
${\cal C}_n$ is a strictly braided tensor category when $n> 2$.

It follows from Theorem \ref {2.2} that

\begin {Corollary} \label {3.5} (Duality Theorem) Let $H$ be a
finite dimensional Hopf algebra with $C_{H,H} = C_{H,H}^{-1}$ in
$({\cal C},C)$. Then

   $$
(R \# H)\# H^{\hat *}   \cong M_n(R)  \hbox { \ \ \  as algebras
in } ({\cal C}, C) $$ in the following three cases:

(i)  $({\cal C}, C)$ is the braided tensor category $({}_B{\cal
M}, C^R)$ determined by the quasi-triangular structure $R$;

(ii) $({\cal C}, C)$  is  the braided tensor category $({}^B{\cal
M}, C^r)$ determined by the coquasi-triangular structure $r$

(iii)  $({\cal C}, C)$ is the braided tensor category ${}^B_B
{\cal YD}$ or  ${}_B {\cal YD} ^B$ of Yetter-Drinfeld modules.
\end {Corollary}

\begin {Example} \label {3.6}
(see \cite  {MR94} )   Let $A$  denote  the anyonic line algebra,
  i.e., $A= {\bf C}\{x \}/\langle x^n \rangle$, where  ${\bf C} \{x\}$
 is a free algebra
  over the complex
 field {\bf C}  and $\langle x^n\rangle$  is an ideal generalized by $x^n$ of
  ${\bf C} \{x\}$.
Set ${\bf C} [\xi]={\bf C} \{x\} / \langle x^n\rangle  $  with
$\xi ^n=0$.
 Its comultiplication, counit,  and antipode are
 $$ \Delta (\xi) = \xi \otimes 1 + 1 \otimes \xi , \hbox { \ \ \  \ \ \ }
 \epsilon (\xi ) =0 , \hbox { \ \ \ and \ \ \ } S(\xi ) =- \xi , $$
 respectively.
 It is straightforward to check that  $H$  is an anyonic Hopf algebra.
 Let $H$ denote the braided group analogue  $\underline {{\bf CZ} _n }$ of   ${\bf CZ} _n $.
Since  $  {\bf CZ} _n $ is commutative, we have  $C_{ {\bf CZ} _n
, {\bf CZ} _n }=C_{ {\bf CZ} _n ,   {\bf CZ} _n }^{-1}$. Let $H$
act on $A$ trivially. By Corollary \ref {3.5}, we have
 $$  ({\bf C} [\xi] \# \underline {{\bf CZ} _n })\#(\underline {{\bf CZ} _n })^{\hat *}
 \cong  M_n({\bf C}[\xi]) \hbox { \ \ \  as algebras in } ({} ^{{\bf CZ} _n}
{\cal M} , C^r).$$

By the way, we also have that the Drinfeld double of  $\underline
{{\bf CZ} _n }$ is  quasi-triangular in
 $ ({} ^{{\bf CZ} _n}
{\cal M} , C^r)$ with the quasi-triangular structure  $[b]$.

 \end {Example}

\vskip 1cm

{\bf Acknowledgement }  The author thanks the referee and the
editor for valuable suggestions.

\newpage
\section {Appendix}
We denote the multiplication, comutiplication, antipode , braiding
and   inverse braiding  by \[
\begin{tangle}
\cu \step[1],\step[2] \cd \step[2], \step[2]\S \step[2],\step[2]\x
\step[2] \hbox { and } \step[2] \xx \step[2], \hbox {respectively
. In particular, }\\
 \hbox {we denote }\step [2] \object {C_{U,V}  }\step
[3]  \hbox { by } \step [2] \X \ ,\step [3] \hbox {for any } \step
[3] \object {U, V = H } \step [3] \hbox { or } \step [1] \object
{H^*}.\end{tangle}
\]

The proof of Lemma2.4\\
Proof.
\[
\begin{tangle}
\object{H\# {H^{\hat{*}}}}\step[6]\object{H\#
{H^{\hat{*}}}}\\
\step\Cu\\
\step\Step\O \lambda\\
\step[3]\object{H\bar{\otimes}{H^{\hat{*}}}}
\end{tangle}
\step=\step
\begin{tangle}
\object{H}\step[2]\object{H^{\hat{*}}}\Step\object{H}\step[3]\object{H^{\hat{*}}}\\
\id\step\cd\step\id\Step\id\Step\coev\\
\id\step\id\Step\X\Step\id\step\cd\step\id\\
\id\step\tu \rightharpoonup\step\cu\step\id\Step\id\step\id\\
\cu\step[3]\XX\Step\id\step\id\\
\step\id\step[3]\ne2\Step\ev\step\id\\
\step\cu\step[7]\id\\
\Step\object{H}\step[8]\object{H^{\hat{*}}}\\
\end{tangle}
\step=\step
\begin{tangle}
\object{H}\step[4]\object{H^{\hat{*}}}\step[4]\object{H}\step[3]\object{H^{\hat{*}}}\\
\id\Step\cd\Step\cd\Step\id\Step\coev\\
\id\Step\id\Step\XX\Step\id\Step\id\step\cd\step\id\\
\id\Step\XX\Step\XX\Step\id\step\id\Step\id\step\id\\
\cu\Step\ev\Step\cu\step\id\Step\id\step\id\\
\step\id\step[8]\XX\Step\id\step\id\\
\step\id\step[7]\ne3\Step\ev\step\id\\
\step\id\step[4]\ne3\step[8]\id\\
\step\cu\step[11]\id\\
\Step\object{H}\step[12]\object{H^{\hat{*}}}\\
\end{tangle} \]
\[
\step=\step
\begin{tangle}
\object{H}\step[4]\object{H^{\hat{*}}}\step[4]\object{H}\step[2]\object{H^{\hat{*}}}\\
\id\Step\cd\Step\cd\step\id\step[3]\Coev\\
\id\Step\id\Step\XX\Step\id\step\id\Step\cd\step[2]\id\\
\id\Step\XX\Step\XX\step\XX\step\cd\step\id\\
\cu\Step\ev\Step\X\Step\X\Step\id\step\id\\
\step\id\step[6]\ne3\step\ev\step\ev\step\id\\
\step\id\step[3]\ne2\step[10]\id\\
\step\cu\step[12]\id\\
\Step\object{H}\step[13]\object{H^{\hat{*}}}\\
\end{tangle} \ \ \ \ \ \ \ \
and \]
\[
\begin{tangle}
\object{H\# {H^{\hat{*}}}}\step[6]\object{H\#
{H^{\hat{*}}}}\\
\O \lambda\step[4]\O \lambda\\
\Cu\\
\step[2]\object{H\bar{\otimes}{H^{\hat{*}}}}
\end{tangle}
\step=\step
\begin{tangle}
\object{H}\Step\object{H^{\hat{*}}}\Step\object{H}\step[3]\object{H^{\hat{*}}}\\
\id\Step\id\Step\id\Step\id\Step\coev\\
\id\Step\id\Step\id\Step\id\step\cd\step\id\\
\id\Step\id\Step\id\Step\X\Step\id\step\id\\
\id\Step\id\Step\cu\step\ev\step\id\\
\id\Step\id\Step\cd\step[4]\id\\
\id\Step\XX\Step\id\step[4]\id\\
\cu\Step\ev\step[4]\id\\
\step\object{H}\step[8]\object{H^{\hat{*}}}\\
\end{tangle}
\step=\step
\begin{tangle}
\object{H}\Step\object{H^{\hat{*}}}\Step\object{H}\step[4]\object{H^{\hat{*}}}\\
\id\Step\id\Step\id\step[3]\id\Step\coev\\
\id\Step\id\Step\id\step[3]\id\step\cd\step\id\\
\id\Step\id\Step\id\step[3]\X\Step\id\step\id\\
\id\Step\id\step\cd\step\cd\ev\step\id\\
\id\Step\id\step\id\Step\X\Step\id\step[3]\id\\
\id\Step\id\step\cu\step\cu\step[3]\id\\
\id\Step\XX\Step\dd\step[4]\id\\
\cu\Step\ev\step[5]\id\\
\step\object{H}\step[10]\object{H^{\hat{*}}}\\
\end{tangle} \]
\[
\step=\step
\begin{tangle}
\object{H}\step[3]\object{H^{\hat{*}}}\Step\object{H}\step[3]\object{H^{\hat{*}}}\\
\id\Step\cd\Step\id\step[3]\id\Step\coev\\
\id\Step\XX\Step\id\step[3]\id\step\cd\step\id\\
\id\Step\id\Step\id\Step\id\step[3]\X\Step\id\step\id\\
\id\Step\id\Step\id\step\cd\step\cd\ev\step\id\\
\id\Step\id\Step\id\step\id\Step\X\Step\id\step[3]\id\\
\id\Step\id\Step\id\step\cu\step\id\Step\id\step[3]\id\\
\id\Step\id\Step\XX\Step\id\Step\id\step[3]\id\\
\id\Step\XX\Step\ev\Step\id\step[3]\id\\
\id\Step\id\Step\nw2\step[4]\ne2\step[3]\id\\
\cu\step[4]\ev\step[5]\id\\
\step\object{H}\step[12]\object{H^{\hat{*}}}\\
\end{tangle}
\step=\step
\begin{tangle}
\object{H}\step[4]\object{H^{\hat{*}}}\step[3]\object{H}\step[3]\object{H^{\hat{*}}}\\
\id\Step\cd\Step\cd\step\id\step[3]\Coev\\
\id\Step\id\Step\XX\Step\id\step\id\Step\cd\step[2]\id\\
\id\Step\XX\Step\XX\step\XX\step\cd\step\id\\
\cu\Step\ev\Step\X\Step\X\Step\id\step\id\\
\step\nw2\step[5]\ne3\step\ev\step\ev\step\id\\
\step[3]\cu\step[10]\id\\
\step[4]\object{H}\step[11]\object{H^{\hat{*}}}\\
\end{tangle}
\step=\step
\begin{tangle}
\object{H\# {H^{\hat{*}}}}\step[6]\object{H\#
{H^{\hat{*}}}}\\
\step\Cu\\
\step\Step\O \lambda\\
\step[3]\object{H\bar{\otimes}{H^{\hat{*}}}}\\
\end{tangle} .
\]
Thus $\lambda $ is an algebraic morphism.

Similarly, we can show that $\rho $ is an anti-algebraic
morphism.\ \ \ \ \begin{picture}(5,5)
\put(0,0){\line(0,1){5}}\put(5,5){\line(0,-1){5}}
\put(0,0){\line(1,0){5}}\put(5,5){\line(-1,0){5}}
\end{picture}\\

The proof of Theorem 2.5.

\[
\begin{tangle}
\object{H\# {H^{\hat{*}}}}\step[6]\object{{H^{\hat{*}}}\#
H}\\
\O \lambda\step[4]\O \rho\\
\Cu\\
\step[2]\object{H\bar{\otimes}{H^{\hat{*}}}}\\
\end{tangle}
\step=\step
\begin{tangle}
\step\object{H}\Step\object{H^{\hat{*}}}\step[3]\object{H^{\hat{*}}}\Step\object{H}\\
\step\id\Step\XX\Step\id\\
\step\XX\Step\XX\\
\cd\step\XX\Step\d\\
\O S\Step\id\step\id\Step\d\step\cd\\
\XX\step\id\step[3]\id\step\XX\\
\id\Step\X\step[3]\X\Step\id\\
\id\Step\id\step\d\step\dd\step\id\Step\id\\
\id\Step\id\Step\X\Step\id\Step\id\\
\id\Step\XX\step\XX\Step\id\\
\d\step\tu \rightharpoonup\step\tu \leftharpoonup\step\dd\\
\step\tu \rho\step[3]\tu \lambda\\
\Step\d\step[4]\id\\
\step[3]\Cu\\
\step[5]\object{H\bar{\otimes}{H^{\hat{*}}}}\\
\end{tangle}\ \ \ .
\Step\Step ......(1)
\]

We show (1) by following five steps. It is easy to check the
following (i) and (ii).

(i)
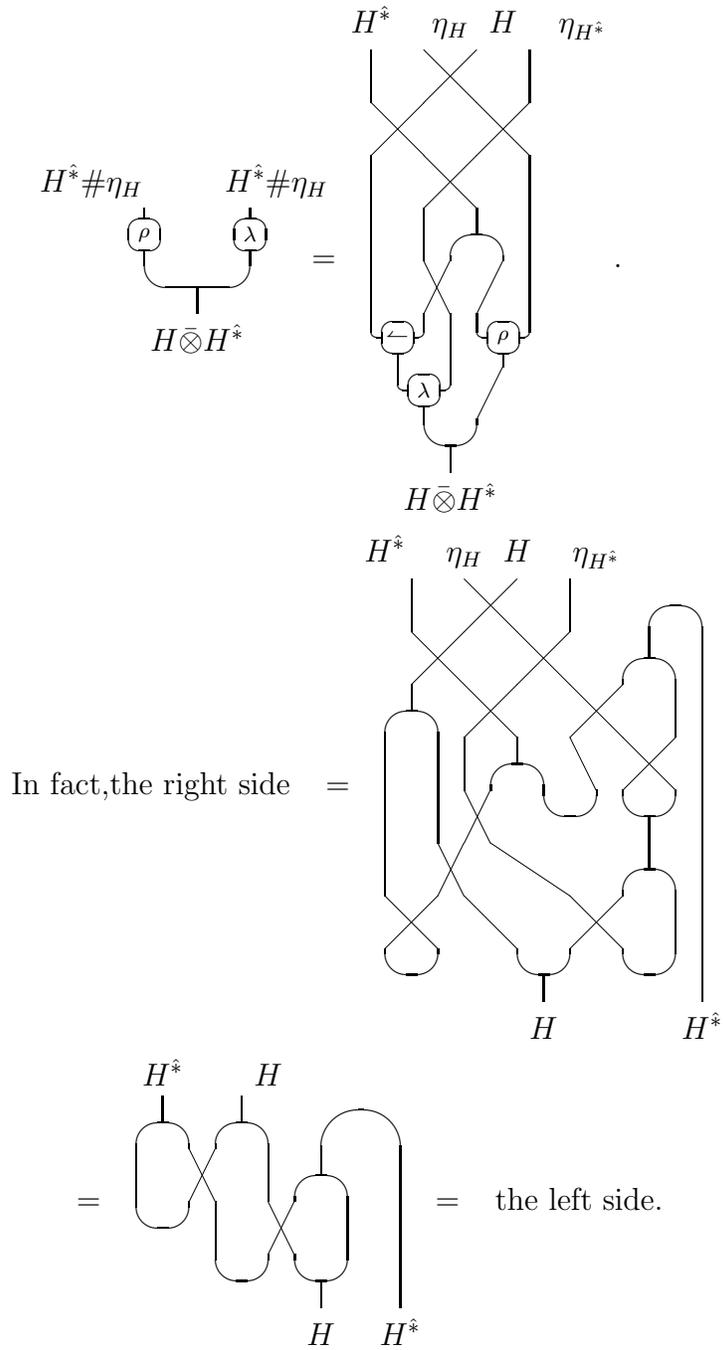
\[
\begin{tangle}
\object{H\#
{\eta_{H^{\hat{*}}}}}\step[7]\object{{\eta_{H^{\hat{*}}}}\#
H}\\
\Step\O \lambda\step[4]\O \rho\\
\Step\Cu\\
\step[4]\object{H\bar{\otimes}{H^{\hat{*}}}}\\
\end{tangle}
\step[3]=\step[3]
\begin{tangle}
\object{H\#
{\eta_{H^{\hat{*}}}}}\step[7]\object{{\eta_{H^{\hat{*}}}}\#
H}\\
\step\d\Step\dd\\
\Step\XX\\
\step\dd\Step\d\\
\step\O \rho\step[4]\O \lambda\\
\step\Cu\\
\step[3]\object{H\bar{\otimes}{H^{\hat{*}}}}\\
\end{tangle}\ \ \ .
\]

(ii)
\[
\begin{tangle}
\object{\eta_H\# {H^{\hat{*}}}}\step[7]\object{H^{\hat{*}}\#
\eta_H}\\
\Step\O \lambda\step[4]\O \rho\\
\Step\Cu\\
\step[4]\object{H\bar{\otimes}{H^{\hat{*}}}}\\
\end{tangle}
\step[3]=\step[3]
\begin{tangle}
\object{\eta_H\# {H^{\hat{*}}}}\step[7]\object{H^{\hat{*}}\#
\eta_H}\\
\step\d\Step\dd\\
\Step\XX\\
\step\dd\Step\d\\
\step\O \rho\step[4]\O \lambda\\
\step\Cu\\
\step[3]\object{H\bar{\otimes}{H^{\hat{*}}}}\\
\end{tangle}\ \ \ .
\]

(iii)
\[
\begin{tangle}
\object{{H^{\hat{*}}}\#\eta_H}\step[7]\object{{H^{\hat{*}}}\#\eta_H
}\\
\Step\O \rho\step[4]\O \lambda\\
\Step\Cu\\
\step[4]\object{H\bar{\otimes}{H^{\hat{*}}}}\\
\end{tangle}
\step=\step
\begin{tangle}
\object{H^{\hat{*}}}\step[3]\object{\eta_H}\Step\object{H}\step[3]\object{\eta_{H^{\hat{*}}}}\\
\id\Step\XX\Step\id\\
\XX\Step\XX\\
\id\Step\XX\Step\id\\
\id\Step\id\step\cd\step\id\\
\id\Step\X\step\dd\step\id\\
\tu \leftharpoonup\step\id\step\tu \rho\\
\step\tu \lambda\step\dd\\
\Step\cu\\
\step[3]\object{H\bar{\otimes}{H^{\hat{*}}}}\\
\end{tangle}\ \ \ .
\]

\[
\hbox {In fact,the right side} \step=\step
\begin{tangle}
\object{H^{\hat{*}}}\step[3]\object{\eta_H}\step[2]\object{H}\step[3]\object{\eta_{H^{\hat{*}}}}\\
\step\id\Step\XX\Step\id\step[3]\coev\\
\step\XX\Step\XX\Step\cd\step\id\\
\cd\step\XX\Step\XX\Step\id\step\id\\
\id\Step\id\step\id\step\cd\step\d\step\XX\step\id\\
\id\Step\id\step\X\Step\ev\step\cu\step\id\\
\id\Step\X\step\nw3\step[4]\cd\step\id\\
\XX\step\nw2\step[3]\XX\Step\id\step\id\\
\ev\step[3]\cu\Step\ev\step\id\\
\step[6]\object{H}\step[6]\object{H^{\hat{*}}}\\
\end{tangle} \]
\[
\step=\step
\begin{tangle}
\step\object{H^{\hat{*}}}\step[4]\object{H}\\
\cd\step\cd\Step\Coev\\
\id\Step\X\Step\id\step\cd\Step\id\\
\ev\step\id\Step\X\Step\id\Step\id\\
\step[3]\ev\step\cu\Step\id\\
\step[7]\object{H}\step[3]\object{H^{\hat{*}}}\\
\end{tangle}
\step=\step \hbox{the left side.}
\]
Thus (iii) holds.

(iv)
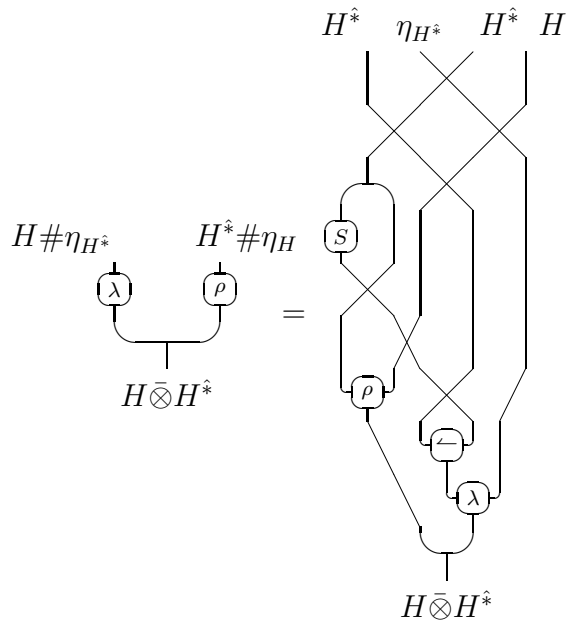
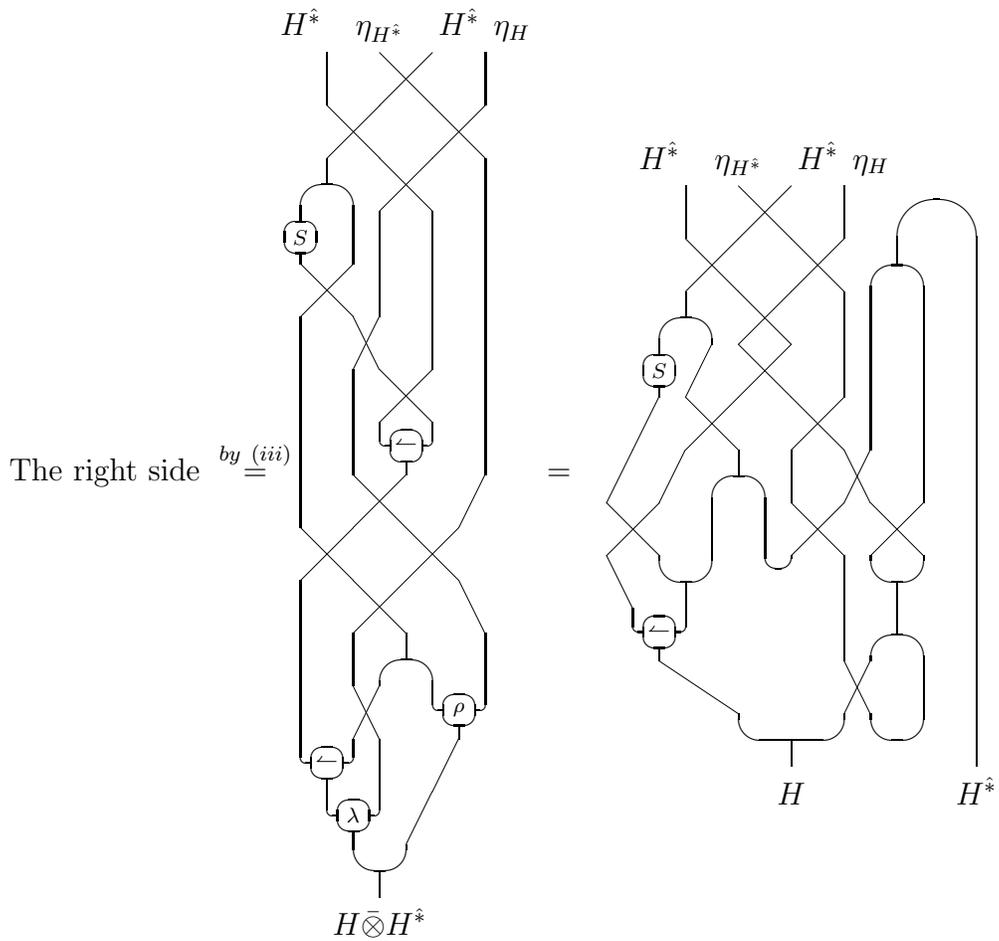
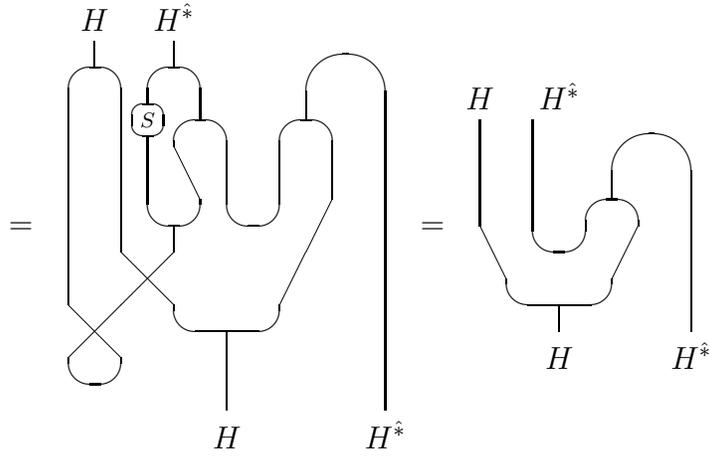
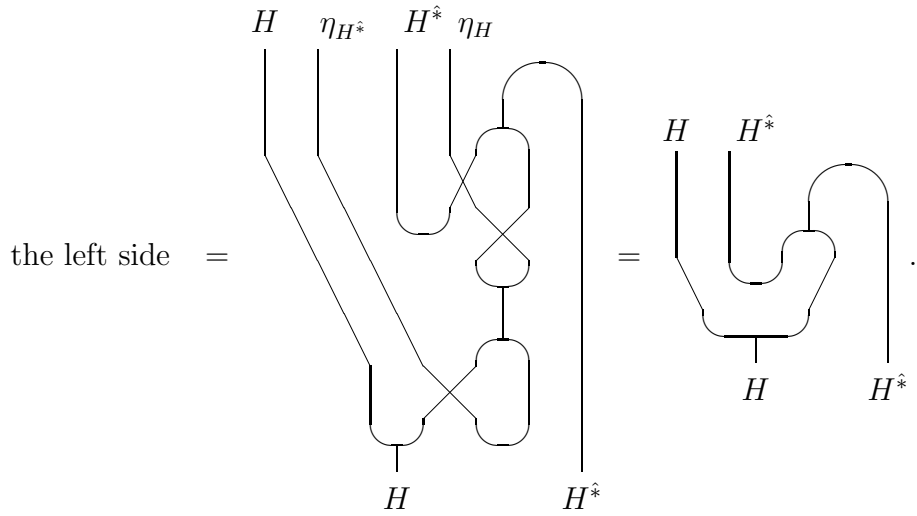
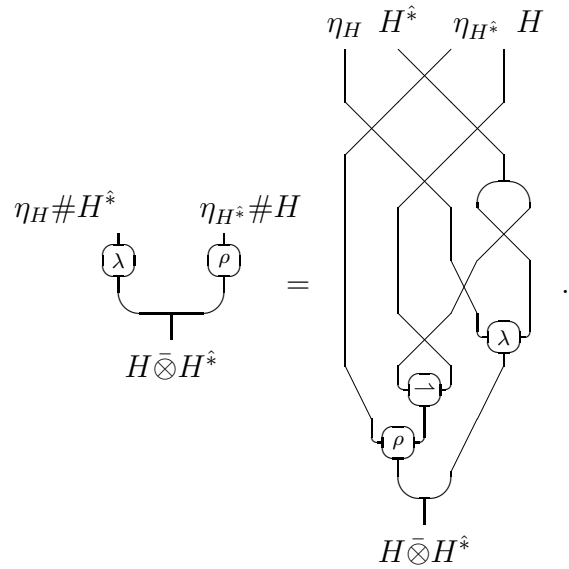
\[
\begin{tangle}
\object{H\# {\eta_{H^{\hat{*}}}}}\step[7]\object{{H^{\hat{*}}}\#
{\eta_H}}\\
\Step\O \lambda\step[4]\O \rho\\
\Step\Cu\\
\step[4]\object{H\bar{\otimes}{H^{\hat{*}}}}\\
\end{tangle}
\step=\step
\begin{tangle}
\object{H^{\hat{*}}}\step[3]\object{\eta_{H^{\hat{*}}}}\step[3]\object{H^{\hat{*}}}\step[2]\object{H}\\
\step\id\Step\XX\Step\id\\
\step\XX\Step\XX\\
\cd\step\XX\Step\id\\
\O S\Step\id\step\id\Step\id\Step\id\\
\XX\step\id\Step\id\Step\id\\
\id\Step\X\Step\id\Step\id\\
\tu \rho\step\XX\step\dd\\
\step\d\step\tu \leftharpoonup\step\id\\
\Step\d\step\tu \lambda\\
\step[3]\cu\\
\step[4]\object{H\bar{\otimes}{H^{\hat{*}}}}\\
\end{tangle}\ \ \ .
\]

\[
\hbox{The right side }\stackrel{by\ (iii)}{=}
\begin{tangle}
\object{H^{\hat{*}}}\step[3]\object{\eta_{H^{\hat{*}}}}\step[3]\object{H^{\hat{*}}}\step[2]\object{\eta_{H}}\\
\step\id\Step\XX\Step\id\\
\step\XX\Step\XX\\
\cd\step\XX\Step\id\\
\O S\Step\id\step\id\Step\id\Step\id\\
\XX\step\id\Step\id\Step\id\\
\id\Step\X\Step\id\Step\id\\
\id\Step\id\step\XX\Step\id\\
\id\Step\id\step\tu \leftharpoonup\Step\id\\
\id\Step\XX\Step\dd\\
\XX\Step\XX\\
\id\Step\XX\Step\d\\
\id\Step\id\step\cd\Step\id\\
\id\Step\X\Step\tu \rho\\
\tu \leftharpoonup\step\id\Step\dd\\
\step\tu \lambda\step\dd\\
\Step\cu\\
\step[3]\object{H\bar{\otimes}{H^{\hat{*}}}}\\
\end{tangle}
\step=\step
\begin{tangle}
\Step\object{H^{\hat{*}}}\step[3]\object{\eta_{H^{\hat{*}}}}\step[3]\object{H^{\hat{*}}}\step[2]\object{\eta_{H}}\\
\step[3]\id\Step\XX\Step\id\Step\Coev\\
\step[3]\XX\Step\XX\step\cd\Step\id\\
\Step\cd\step\XX\Step\id\step\id\Step\id\Step\id\\
\Step\O S\step\dd\step\XX\Step\id\step\id\Step\id\Step\id\\
\step\dd\step\XX\Step\XX\step\id\Step\id\Step\id\\
\dd\step\dd\step\cd\step\id\Step\X\Step\id\Step\id\\
\XX\Step\id\Step\id\step\XX\step\XX\Step\id\\
\d\step\cu\Step\hev\Step\id\step\cu\Step\id\\
\step\tu \leftharpoonup\step[6]\id\step\cd\Step\id\\
\Step\nw3\step[6]\X\Step\id\Step\id\\
\step[5]\Cu\step\ev\Step\id\\
\step[7]\object{H}\step[7]\object{H^{\hat{*}}}\\
\end{tangle}\]
\[
\step=\step
\begin{tangle}
\step\object{H}\step[3]\object{H^{\hat{*}}}\\
\cd\step\cd\step[4]\Coev\\
\id\Step\id\step\O S\step\cd\Step\cd\Step\id\\
\id\Step\id\step\id\step\d\step\id\Step\id\Step\id\Step\id\\
\id\Step\id\step\cu\step\ev\step\dd\Step\id\\
\id\Step\XX\step[4]\dd\step[3]\id\\
\XX\Step\Cu\step[4]\id\\
\ev\step[4]\id\step[6]\id\\
\step[6]\object{H}\step[6]\object{H^{\hat{*}}}\\
\end{tangle}
\step=\step
\begin{tangle}
\object{H}\step[3]\object{H^{\hat{*}}}\\
\id\Step\id\step[3]\Coev\\
\id\Step\id\Step\cd\Step\id\\
\d\step\ev\step\dd\Step\id\\
\step\Cu\step[3]\id\\
\step[3]\object{H}\step[5]\object{H^{\hat{*}}}\\
\end{tangle}
\]

and
\[
\hbox{the left side}\step=\step
\begin{tangle}
\object{H}\step[3]\object{\eta_{H^{\hat{*}}}}\step[3]\object{H^{\hat{*}}}\step[2]\object{\eta_H}\\
\id\Step\id\step[3]\id\Step\id\Step\Coev\\
\id\Step\id\step[3]\id\Step\id\step\cd\Step\id\\
\d\step\d\step[2]\id\Step\X\Step\id\Step\id\\
\step\d\step\d\step\ev\step\XX\Step\id\\
\Step\d\step\d\step[3]\cu\Step\id\\
\step[3]\d\step\d\Step\cd\Step\id\\
\step[4]\id\Step\XX\Step\id\Step\id\\
\step[4]\cu\Step\ev\Step\id\\
\step[5]\object{H}\step[7]\object{H^{\hat{*}}}\\
\end{tangle}
\step=\step
\begin{tangle}
\object{H}\step[3]\object{H^{\hat{*}}}\\
\id\Step\id\step[3]\Coev\\
\id\Step\id\Step\cd\Step\id\\
\d\step\ev\step\dd\Step\id\\
\step\Cu\step[3]\id\\
\step[3]\object{H}\step[5]\object{H^{\hat{*}}}\\
\end{tangle}\ \ .
\]

Thus (iv) holds.

(v)
\[
\begin{tangle}
\object{\eta_H\#
{H^{\hat{*}}}}\step[7]\object{\eta_{H^{\hat{*}}}\#
H}\\
\Step\O \lambda\step[4]\O \rho\\
\Step\Cu\\
\step[4]\object{H\bar{\otimes}{H^{\hat{*}}}}\\
\end{tangle}
\step=\step
\begin{tangle}
\object{\eta_H}\Step\object{H^{\hat{*}}}\step[3]\object{\eta_{H^{\hat{*}}}}\step[2]\object{H}\\
\id\Step\XX\Step\id\\
\XX\Step\XX\\
\id\Step\XX\step\cd\\
\id\Step\id\Step\id\step\XX\\
\id\Step\id\Step\X\Step\id\\
\id\Step\XX\step\tu \lambda\\
\d\step\tu \rightharpoonup\step\dd\\
\step\tu \rho\step\dd\\
\Step\cu\\
\step[3]\object{H\bar{\otimes}{H^{\hat{*}}}}\\
\end{tangle}\ \ \ .
\]

\[
\hbox{In fact, the right side}\step=\step
\begin{tangle}
\object{\eta_H}\Step\object{H^{\hat{*}}}\step[3]\object{\eta_{H^{\hat{*}}}}\step[2]\object{H}\\
\id\Step\XX\Step\id\step[3]\Coev\\
\XX\Step\XX\Step\cd\Step\id\\
\id\Step\XX\step\cd\step\id\Step\id\Step\id\\
\id\step\cd\step\id\step\XX\step\id\Step\id\Step\id\\
\id\step\id\Step\id\step\X\Step\X\Step\id\Step\id\\
\id\step\d\step\X\step\cu\step\ev\Step\id\\
\d\step\d\hev\step\cd\step[5]\id\\
\step\id\Step\XX\Step\id\step[5]\id\\
\step\ev\Step\XX\step[5]\id\\
\step[5]\cu\step[5]\id\\
\step[6]\object{H}\step[6]\object{H^{\hat{*}}}\\
\end{tangle}
\step=\step
\begin{tangle}
\step\object{H^{\hat{*}}}\step[3]\object{H}\\
\cd\step\cd\Step\Coev\\
\XX\step\XX\step\cd\Step\id\\
\id\Step\X\Step\id\step\id\Step\id\Step\id\\
\ev\step\XX\step\id\Step\id\Step\id\\
\step[3]\id\Step\X\Step\id\Step\id\\
\step[3]\XX\step\ev\step[2]\id\\
\step[3]\cu\step[5]\id\\
\step[4]\object{H}\step[6]\object{H^{\hat{*}}}\\
\end{tangle}
\]

and
\[
\hbox{the left side}\step=\step
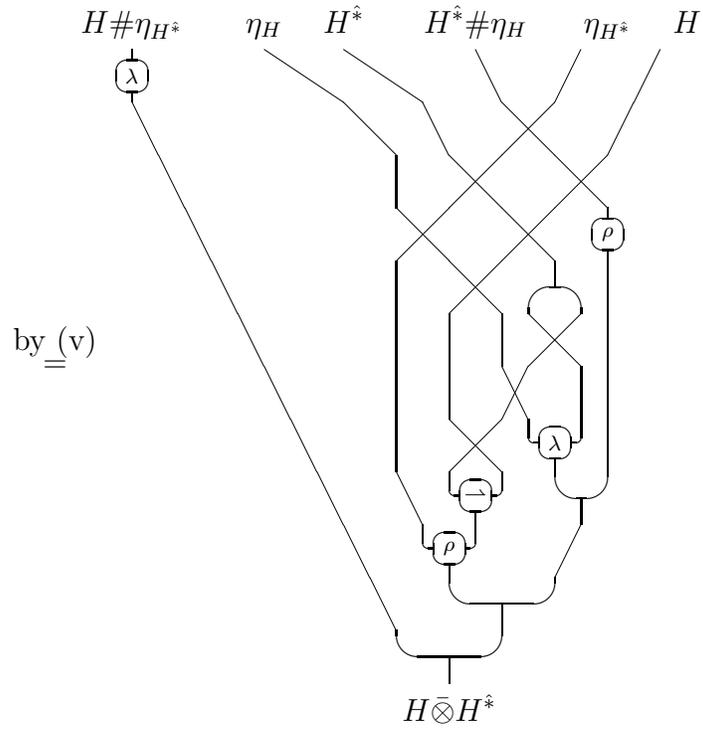
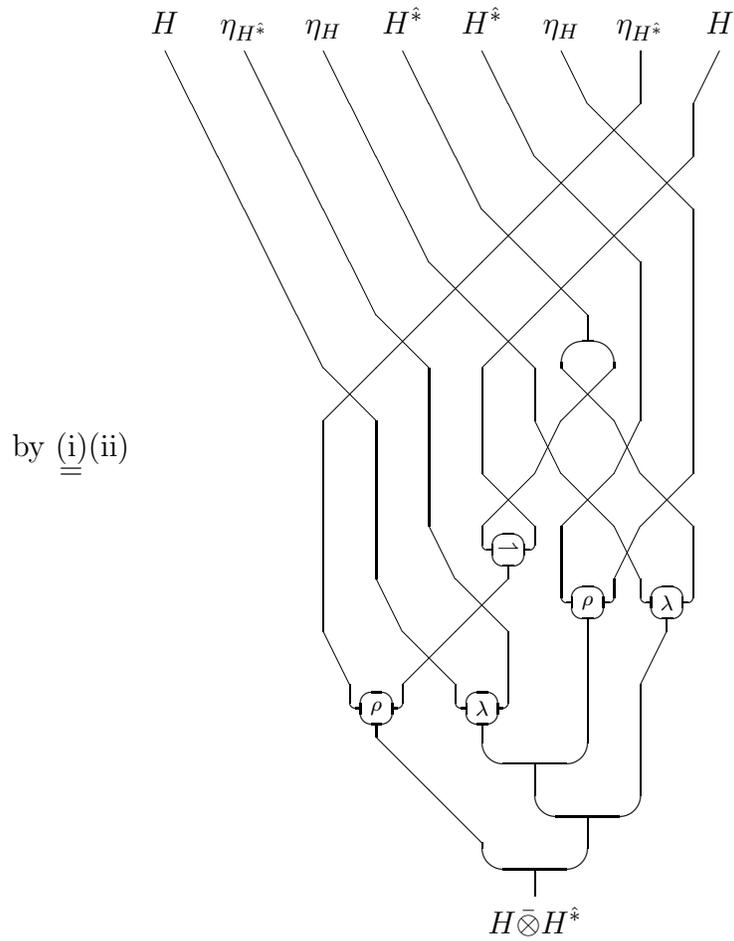
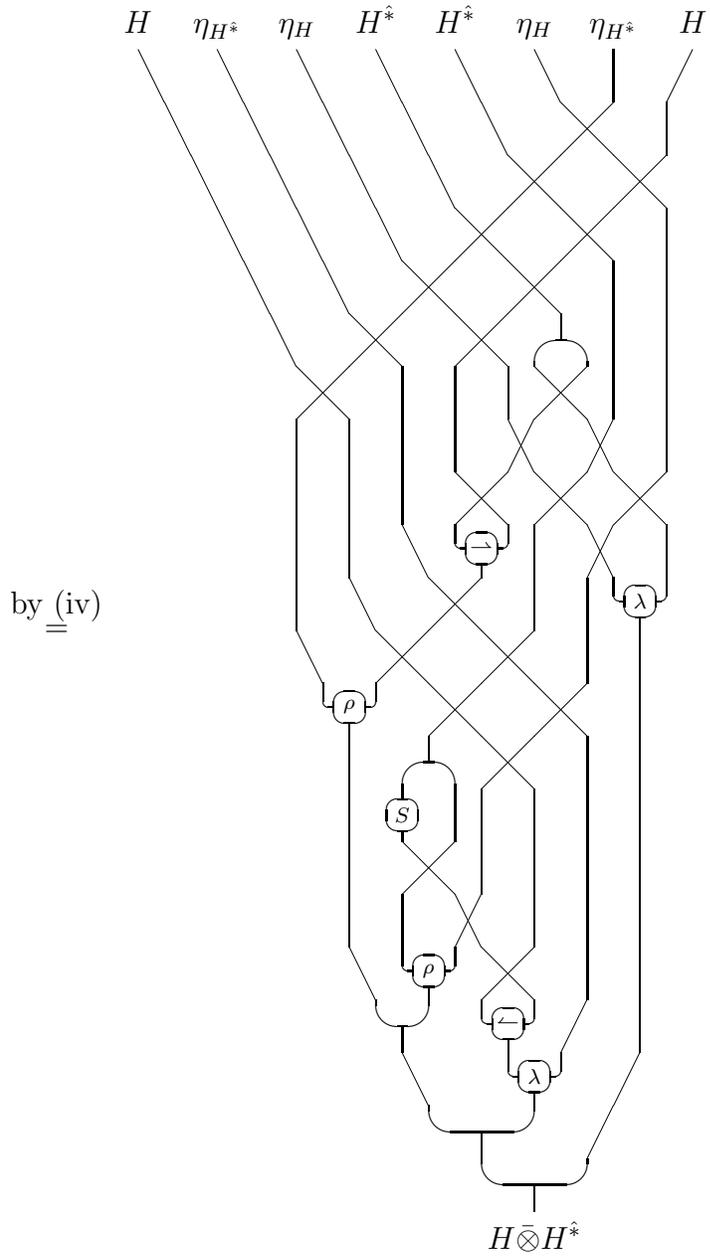
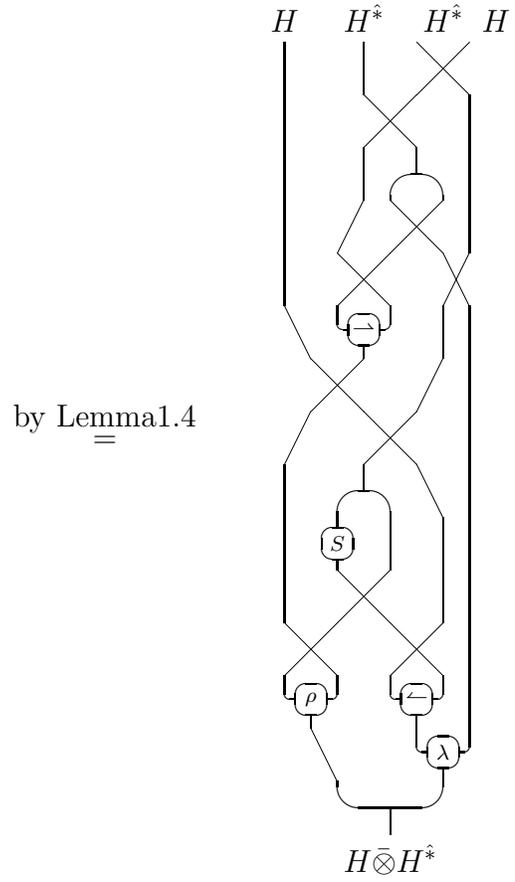
\begin{tangle}
\object{H^{\hat{*}}}\step[3]\object{H}\\
\id\Step\id\Step\coev\\
\id\Step\XX\Step\id\\
\id\Step\cu\Step\id\\
\id\Step\cd\Step\id\\
\XX\Step\id\Step\id\\
\id\Step\ev\Step\id\\
\object{H}\step[6]\object{H^{\hat{*}}}\\
\end{tangle}
\step=\step
\begin{tangle}
\object{H^{\hat{*}}}\step[3]\object{H}\\
\id\Step\id\Step\Coev\\
\id\Step\XX\step[3]\id\\
\id\Step\id\Step\d\step[2]\id\\
\id\step\cd\step\cd\step\id\\
\id\step\id\Step\X\Step\id\step\id\\
\id\step\cu\step\cu\step\id\\
\XX\Step\dd\step[2]\id\\
\id\Step\ev\step[3]\id\\
\object{H}\step[7]\object{H^{\hat{*}}}\\
\end{tangle}\]
\[
\step=\step
\begin{tangle}
\step\object{H^{\hat{*}}}\step[3]\object{H}\\
\cd\step[2]\id\Step\Coev\\
\id\Step\id\step[2]\XX\step[3]\id\\
\id\step[2]\id\Step\id\Step\d\Step\id\\
\id\step[2]\id\step\cd\step\cd\step\id\\
\id\Step\id\step\id\Step\X\Step\id\step\id\\
\id\Step\id\step\cu\step\id\Step\id\step\id\\
\id\Step\XX\Step\id\Step\id\step\id\\
\XX\Step\XX\Step\id\step\id\\
\id\Step\ev\Step\ev\step\id\\
\object{H}\step[9]\object{H^{\hat{*}}}\\
\end{tangle}
\step=\step
\begin{tangle}
\step\object{H^{\hat{*}}}\step[3]\object{H}\\
\cd\step\cd\Step\Coev\\
\id\Step\id\step\id\Step\id\step\cd\Step\id\\
\id\Step\id\step\id\Step\X\Step\id\Step\id\\
\id\Step\id\step\XX\step\id\Step\id\Step\id\\
\id\Step\id\step\cu\step\id\Step\id\Step\id\\
\id\Step\XX\Step\id\step\dd\Step\id\\
\XX\Step\ev\dd\step[3]\id\\
\id\Step\d\step[3]\id\step[4]\id\\
\id\step[3]\Ev\step[6]\id\\
\object{H}\step[10]\object{H^{\hat{*}}}\\
\end{tangle}
\step=\step\hbox{the right side.}
\]
Thus (v) holds.

Now we show that the relation (1) holds.
\[
\hbox{the left side of (1)}\stackrel{\hbox{by
Lemma1.4}}{\step=\step[5]}
\begin{tangle}
\object{H\# {\eta_{H^{\hat{*}}}}}\step[8]\object{{\eta_H}\#
H^{\hat{*}}}\step[9]\object{H^{\hat{*}}\#\eta_H}\step[6]\object{\eta_{H^{\hat{*}}}\# H}\\
\O \lambda\step[6]\Step\O \lambda\step[11]\XX \\
\nw3\step[7]\nw3\step[10]\O \rho\step[2]\O \rho\\
\step[3]\nw3\step[7]\nw3\step[7]\cu\\
\step[6]\nw3\step[7]\d\step[4]\dd\\
\step[9]\nw3\step[5]\Cu\\
\step[12]\d\step[4]\id\\
\step[13]\Cu\\
\step[15]\object{H\bar{\otimes}{H^{\hat{*}}}}\\
\end{tangle} \]
\[
\step\stackrel{\hbox{by (v)}}{=}\step
\begin{tangle}
\object{H\#
{\eta_{H^{\hat{*}}}}}\step[5]\object{\eta_H}\step[3]\object{H^{\hat{*}}}\step[5]
\object{H^{\hat{*}}\#\eta_H}\step[5]\object{\eta_{H^{\hat{*}}}}\step[3]\object{H}\\
\O \lambda\step[5]\nw3\step[2]\nw3\step[4]\d\Step\dd\Step\dd\\
\d\step[7]\nw2\Step\d\Step\XX\Step\dd\\
\step\d\step[8]\id\Step\XX\Step\XX\\
\Step\d\step[7]\XX\Step\XX\Step\O \rho\\
\step[3]\d\step[6]\id\Step\XX\step\cd\step\id\\
\step[4]\d\step[5]\id\Step\id\Step\id\step\XX\step\id\\
\step[5]\d\step[4]\id\Step\id\Step\X\Step\id\step\id\\
\step[6]\d\step[3]\id\Step\XX\step\tu \lambda\step\id\\
\step[7]\d\Step\d\step\tu \rightharpoonup\Step\cu\\
\step[8]\d\Step\tu \rho\step[3]\dd\\
\step[9]\d\Step\Cu\\
\step[10]\Cu\\
\step[12]\object{H\bar{\otimes}{H^{\hat{*}}}}\\
\end{tangle}
\]

\[
\step\stackrel{\hbox{by (i)(ii)}}{=}\step
\begin{tangle}
\object{H}\step[3]\object{\eta_{H^{\hat{*}}}}\step[3]\object{\eta_H}\step[3]\object{H^{\hat{*}}}
\step[3]\object{H^{\hat{*}}}\step[3]\object{\eta_H}\step[3]\object{\eta_{H^{\hat{*}}}}\step[3]\object{H}\\
\d\Step\d\Step\d\Step\d\Step\d\Step\d\Step\id\Step\dd\\
\step\d\Step\d\Step\d\Step\d\Step\d\Step\XX\Step\id\\
\Step\d\Step\d\Step\d\Step\d\Step\XX\Step\XX\\
\step[3]\d\Step\d\Step\d\Step\XX\Step\XX\Step\id\\
\step[4]\d\Step\d\Step\XX\Step\XX\Step\id\Step\id\\
\step[5]\d\Step\XX\Step\XX\step\cd\step\id\Step\id\\
\step[6]\XX\Step\id\Step\id\Step\id\step\XX\step\id\Step\id\\
\step[6]\id\Step\id\Step\id\Step\id\Step\X\Step\X\Step\id\\
\step[6]\id\Step\id\Step\id\Step\XX\step\XX\step\XX\\
\step[6]\id\Step\id\Step\d\step\tu
\rightharpoonup\step\id\Step\X\Step\id\\
\step[6]\id\Step\d\Step\XX\Step\tu \rho\step\tu \lambda\\
\step[6]\d\Step\XX\Step\id\step[3]\id\Step\dd\\
\step[7]\tu \rho\Step\tu \lambda\step[3]\id\Step\id\\
\step[8]\nw2\step[3]\Cu\Step\id\\
\step[10]\nw2\step[3]\Cu\\
\step[12]\Cu\\
\step[14]\object{H\bar{\otimes}{H^{\hat{*}}}}\\
\end{tangle}
\]
\[
\step\stackrel{\hbox{by (iv)}}{=}\step
\begin{tangle}
\object{H}\step[3]\object{\eta_{H^{\hat{*}}}}\step[3]\object{\eta_H}\step[3]\object{H^{\hat{*}}}
\step[3]\object{H^{\hat{*}}}\step[3]\object{\eta_H}\step[3]\object{\eta_{H^{\hat{*}}}}\step[3]\object{H}\\
\d\Step\d\Step\d\Step\d\Step\d\Step\d\Step\id\Step\dd\\
\step\d\Step\d\Step\d\Step\d\Step\d\Step\XX\Step\id\\
\Step\d\Step\d\Step\d\Step\d\Step\XX\Step\XX\\
\step[3]\d\Step\d\Step\d\Step\XX\Step\XX\Step\id\\
\step[4]\d\Step\d\Step\XX\Step\XX\Step\id\Step\id\\
\step[5]\d\Step\XX\Step\XX\step\cd\step\id\Step\id\\
\step[6]\XX\Step\id\Step\id\Step\id\step\XX\step\id\Step\id\\
\step[6]\id\Step\id\Step\id\Step\id\Step\X\Step\X\Step\id\\
\step[6]\id\Step\id\Step\id\Step\XX\step\XX\step\XX\\
\step[6]\id\Step\id\Step\d\step\tu
\rightharpoonup\step\id\Step\X\Step\id\\
\step[6]\id\Step\d\Step\XX\Step\id\Step\id\step\tu \lambda\\
\step[6]\d\Step\XX\Step\XX\Step\id\Step\id\\
\step[7]\tu \rho\Step\XX\Step\XX\Step\id\\
\step[8]\id\Step\cd\step\XX\Step\id\Step\id\\
\step[8]\id\Step\O S\Step\id\step\id\Step\id\Step\id\Step\id\\
\step[8]\id\Step\XX\step\id\Step\id\Step\id\Step\id\\
\step[8]\id\Step\id\Step\X\Step\id\Step\id\Step\id\\
\step[8]\d\step\tu \rho\step\XX\Step\id\Step\id\\
\step[9]\cu\Step\tu \leftharpoonup\step\dd\Step\id\\
\step[10]\d\step[3]\tu \lambda\step[2]\dd\\
\step[11]\Cu\step[2]\dd\\
\step[13]\Cu\\
\step[15]\object{H\bar{\otimes}{H^{\hat{*}}}}\\
\end{tangle}
\step\stackrel{\hbox{by Lemma1.4}}{=}\step[3]
\begin{tangle}
\object{H}\step[3]\object{H^{\hat{*}}}\step[3]\object{H^{\hat{*}}}\Step\object{H}\\
\id\step[3]\id\Step\XX\\
\id\step[3]\XX\Step\id\\
\id\step[3]\id\step\cd\step\id\\
\id\Step\dd\step\XX\step\id\\
\id\Step\XX\Step\X\\
\d\step\tu \rightharpoonup\Step\id\step\id\\
\step\XX\Step\dd\step\id\\
\dd\Step\XX\Step\id\\
\id\Step\cd\step\d\step\id\\
\id\Step\O S\Step\id\Step\id\step\id\\
\id\Step\XX\Step\id\step\id\\
\XX\Step\XX\step\id\\
\tu \rho\Step\tu \leftharpoonup\step\id\\
\step\d\step[3]\tu \lambda\\
\Step\Cu\\
\step[4]\object{H\bar{\otimes}{H^{\hat{*}}}}\\
\end{tangle}
\]
\[
\step=\step
\begin{tangle}
\step\object{H}\step[2]\object{H^{\hat{*}}}\step[3]\object{H^{\hat{*}}}\Step\object{H}\\
\step\id\Step\XX\Step\id\\
\step\XX\Step\XX\\
\cd\step\XX\step\cd\\
\O S\Step\id\step\id\Step\id\step\XX\\
\XX\step\id\Step\X\Step\id\\
\id\Step\X\step\dd\step\id\Step\id\\
\id\step\dd\step\X\Step\id\Step\id\\
\id\step\XX\step\XX\Step\id\\
\id\step\tu \rightharpoonup\step\tu \leftharpoonup\step\dd\\
\tu \rho\step[3]\tu \lambda\\
\step\id\step[4]\dd\\
\step\Cu\\
\step[3]\object{H\bar{\otimes}{H^{\hat{*}}}}\\
\end{tangle}
\step=\step
\begin{tangle}
\hbox{the right side of (1). }
\end{tangle}
\ \ \ \
\begin{picture}(5,5)
\put(0,0){\line(0,1){5}}\put(5,5){\line(0,-1){5}}
\put(0,0){\line(1,0){5}}\put(5,5){\line(-1,0){5}}
\end{picture}\\
\]

The proof of Theorem 2.8 .\\
(i)
\[
\step Let\step[3]
\begin{tangle}
\object{H^{\hat{*}}}\\
\id\\
\O w\\
\id\\
\object{H\# H^{\hat{*}}}\\
\end{tangle}
\step=\step
\begin{tangle}
\object{H^{\hat{*}}}\step[3]\object{\eta_H}\\
\step[0.5]\O { \bar S}\Step\id\\
\step[0.5]\tu \rho\\
\obox 3{\lambda^{-1}}\\
\step[1.5]\id\\
\step[1.5]\object{H\# H^{\hat{*}}}\\
\end{tangle}
\Step and\Step
\begin{tangle}
\step\object{R}\\
\td \phi\\
\object{H^{\hat{*}}}\Step\object{R}\\
\end{tangle}
\step=\step
\begin{tangle}
\step[4]\object{R}\\
\ro {b'} \Step\id\\
\id\Step\tu \alpha\\
\object{H^{\hat{*}}}\step[3]\object{R}\\
\end{tangle}
\step=\step
\begin{tangle}
\step[4]\object{R}\\
\coev\Step\id\\
\XX\Step\id\\
\id\Step\tu \alpha\\
\object{H^{\hat{*}}}\step[3]\object{R}\\
\end{tangle}\ \ \ , \hbox {where } \step [4] \object { b' = C_{H,
H^*}b_H .}
\]

We define

\[
\begin{tangle}
\object{R\otimes(H\# H^{\hat{*}})}\\
\id\\
\O \Psi\\
\id\\
\object{(R\# H)\# H^{\hat{*}}}\\
\end{tangle}
\step[3]=\step[3]
\begin{tangle}
\step\object{R}\step[4]\object{H}\Step\object{H^{\hat{*}}}\\
\td \phi\step[3]\id\Step\id\\
\O S\Step\id\step[3]\id\Step\id\\
\x\step[3]\id\Step\id\\
\id\step\td w\Step\id\Step\id\\
\id\step\id\step\cd\step\id\Step\id\\
\id\step\id\step\id\Step\X\Step\id\\
\id\step\id\step\tu \rightharpoonup\step\cu\\
\id\step\cu\step[3]\id\\
\object{R}\Step\object{H}\step[4]\object{H^{\hat{*}}}\\
\end{tangle}
\step[3]and\step[5]
\begin{tangle}
\object{(R\# H)\# H^{\hat{*}}}\\
\id\\
\O \Phi\\
\id\\
\object{R\otimes(H\# H^{\hat{*}})}\\
\end{tangle}
\step[3]=\step[3]
\begin{tangle}
\step\object{R}\step[4]\object{H}\Step\object{H^{\hat{*}}}\\
\td \phi\step[3]\id\Step\id\\
\x\step[3]\id\Step\id\\
\id\step\td w\Step\id\Step\id\\
\id\step\id\step\cd\step\id\Step\id\\
\id\step\id\step\id\Step\X\Step\id\\
\id\step\id\step\tu \rightharpoonup\step\cu\\
\id\step\cu\step[3]\id\\
\object{R}\Step\object{H}\step[4]\object{H^{\hat{*}}}\\
\end{tangle}\ \ \ .
\]

We see that
\[
\Psi\Phi\step=\step
\begin{tangle}
\step[7]\object{R}\step[5]\object{H\# H^{\hat{*}}}\\
\ro {b'} \step\ro {b'}\Step\id\step[4]\dd\\
\O S\Step\id\step\d\step\tu \alpha\step[3]\dd\\
\O w\Step\id\Step\x\step[3]\dd\\
\d\step\tu \alpha\Step\O w\Step\dd\\
\step\x\Step\dd\step\dd\\
\step\id\Step\cu\step\dd\\
\step\id\step[3]\cu\\
\step\object{R}\step[4]\object{H\# H^{\hat{*}}}\\
\end{tangle}
\step \stackrel{\hbox{ since } w \hbox{ is algebraic }}{=}\step
\begin{tangle}
\step[7]\object{R}\step[4]\object{H\# H^{\hat{*}}}\\
\step[3]\coev\Step\id\step[4]\id\\
\coev\step\d\step\tu \alpha\step[4]\id\\
\O S\Step\id\Step\x\step[4]\dd\\
\d\step\tu \alpha\Step\id\step[3]\dd\\
\step\x\Step\dd\step[2]\dd\\
\step\id\Step\cu\step[2]\dd\\
\step\id\step[3]\O w\Step\dd\\
\step\id\step[3]\cu\\
\step\object{R}\step[5]\object{H\# H^{\hat{*}}}\\
\end{tangle} \]
\[
\step=\step
\begin{tangle}
\step[6]\object{R}\step[4]\object{H\# H^{\hat{*}}}\\
\coev\step\coev\step\id\step[4]\id\\
\O S\Step\X\Step\id\step\id\step[4]\id\\
\cu\step\cu\step\id\step[3]\dd\\
\step\O w\step[3]\tu \alpha\Step\dd\\
\step\d\Step\dd\Step\dd\\
\Step\x\Step\dd\\
\Step\id\step[2]\cu\\
\Step\object{R}\step[4]\object{H\# H^{\hat{*}}}\\
\end{tangle}
\step=\step
\begin{tangle}
\step[5]\object{R}\step[4]\object{H\# H^{\hat{*}}}\\
\step\coev\Step\id\step[4]\id\\
\cd\step\tu \alpha\step[4]\id\\
\O S\Step\id\Step\id\step[5]\id\\
\cu\Step\id\step[5]\id\\
\step\O w\Step\dd\step[4]\dd\\
\step\x\step[4]\dd\\
\step\id\Step\Cu\\
\step\object{R}\step[4]\object{H\# H^{\hat{*}}}\\
\end{tangle}
\step=\step
\begin{tangle}
\object{R}\step[4]\object{H\# H^{\hat{*}}}\\
\id\step[4]\id\\
\id\step[4]\id\\
\object{R}\step[4]\object{H\# H^{\hat{*}}}\\
\end{tangle} \ \ \ \ \ .
\]
Similarly, we have  $\Phi \Psi = id $. Thus $\Phi$ is invertible.

Now we show that $\Phi$ is algebraic.
\[
Let\step[12]
\begin{tangle}
\step\object{(R\# H)\# H^{\hat{*}}}\\
\step\id\\
\obox 2{\Phi'}\\
\step\id\\
\step\object{R\otimes(H\bar{\otimes}H^{\hat{*}})}\\
\end{tangle}
\step[3]=\step[3]
\begin{tangle}
\step[4]\object{R}\step[3]\object{\eta_H}\step[4]\object{H\# H^{\hat{*}}}\\
\ro {b'}\Step\id\step[3]\id\step[4]\id\\
\O {\bar{S}}\Step\tu \alpha\step[3]\id\step[4]\O \lambda\\
\d\Step\id\step[3]\dd\step[3]\dd\\
\step\x\Step\dd\step[3]\dd\\
\step\id\Step\tu \rho\step[3]\dd\\
\step\id\step[3]\Cu\\
\step\object{R}\step[5]\object{H\bar{\otimes}H^{\hat{*}}}\\
\end{tangle}\ \ \ .
\]

It is clear that $\Phi = (id \otimes \lambda ^{-1})\Phi '$.
Consequently, we only need show that $\Phi '$ is alegebraic.

\[
Let\step[12]
\begin{tangle}
\object{R}\\
\id\\
\O \xi\\
\id\\
\object{R\otimes(H\bar{\otimes}H^{\hat{*}})}\\
\end{tangle}
\step[3]=\step[3]
\begin{tangle}
\step[4]\object{R}\step[3]\object{\eta_H}\\
\ro {b'}\Step\id\step[3]\id\\
\O {\bar{S}}\Step\tu \alpha\Step\dd\\
\d\Step\id\Step\dd\\
\step\x\step[2]\id\\
\step\id\Step\tu \rho\\
\step\object{R}\step[4]\object{H\bar{\otimes}H^{\hat{*}}}\\
\end{tangle}\ \ \ .
\]

We have that

\[
\begin{tangle}
\step\object{(R\# H)\# H^{\hat{*}}}\\
\step\id\\
\obox 2{\Phi'}\\
\step\id\\
\step\object{R\otimes(H\bar{\otimes}H^{\hat{*}})}\\
\end{tangle}
\step=\step[3]
\begin{tangle}
\step\object{R}\step[5]\object{H\# H^{\hat{*}}}\\
\td \xi\step[4]\O \lambda\\
\id\Step\Cu\\
\object{R}\step[5]\object{H\bar{\otimes} H^{\hat{*}}}\\
\end{tangle}\ \ \ .
\]

We claim that

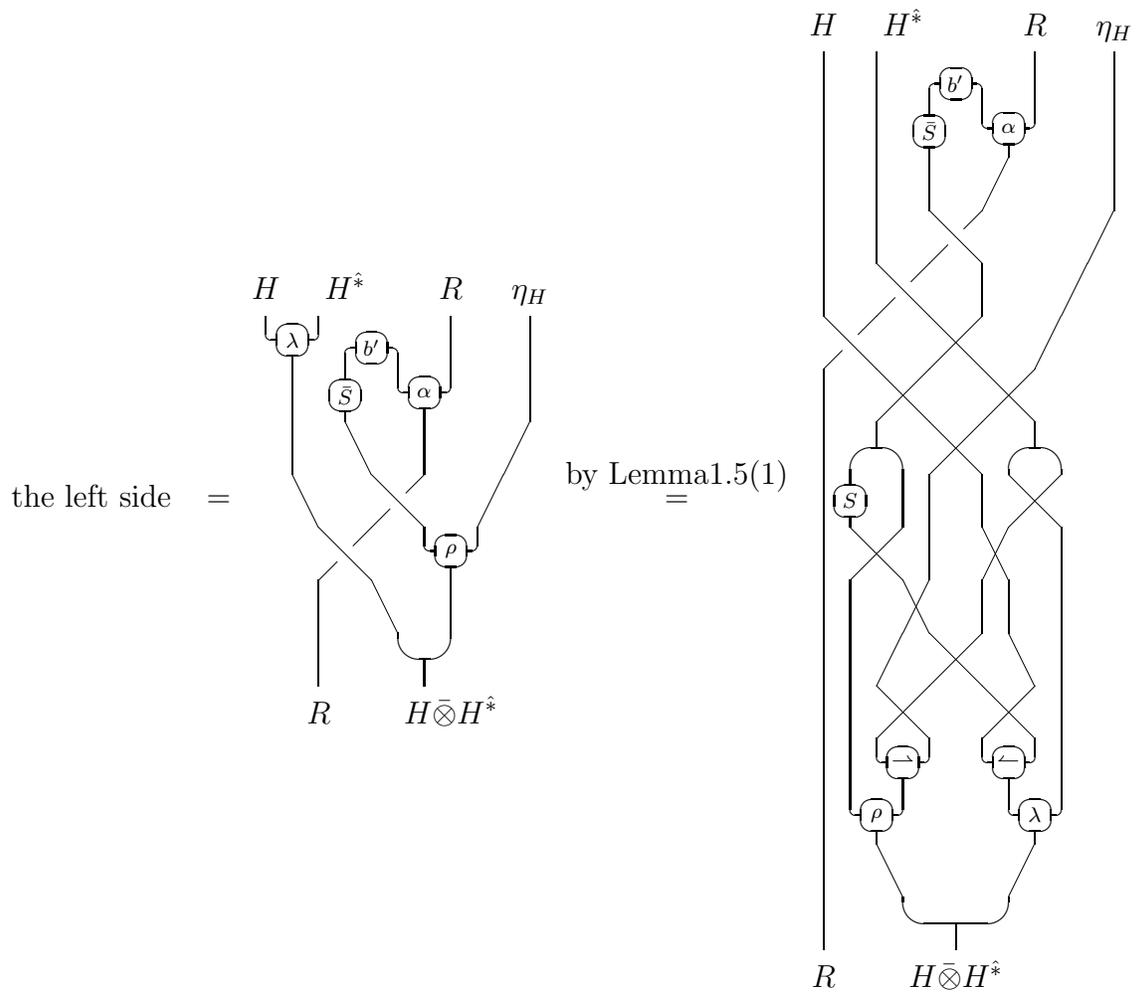
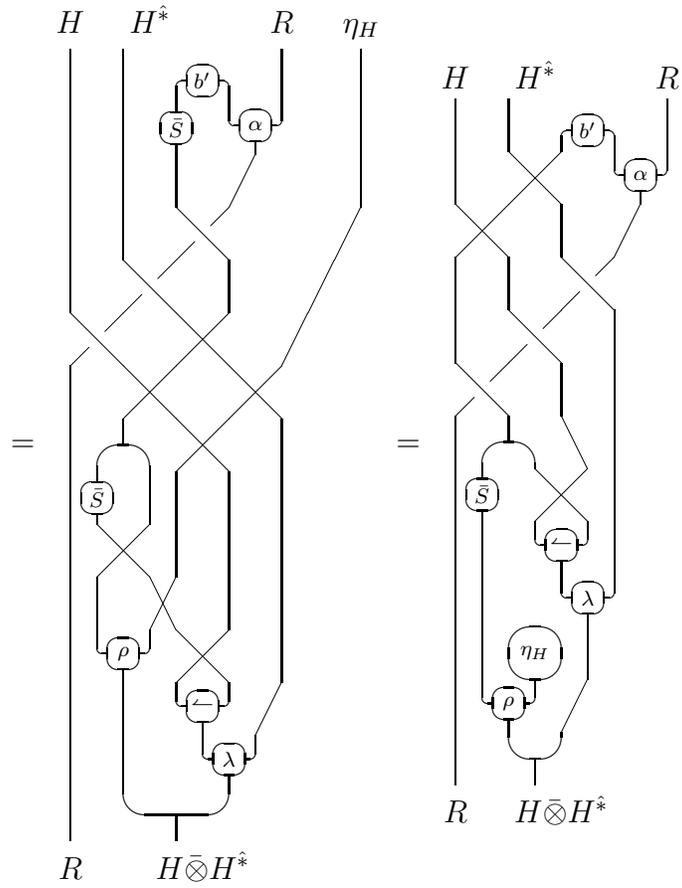
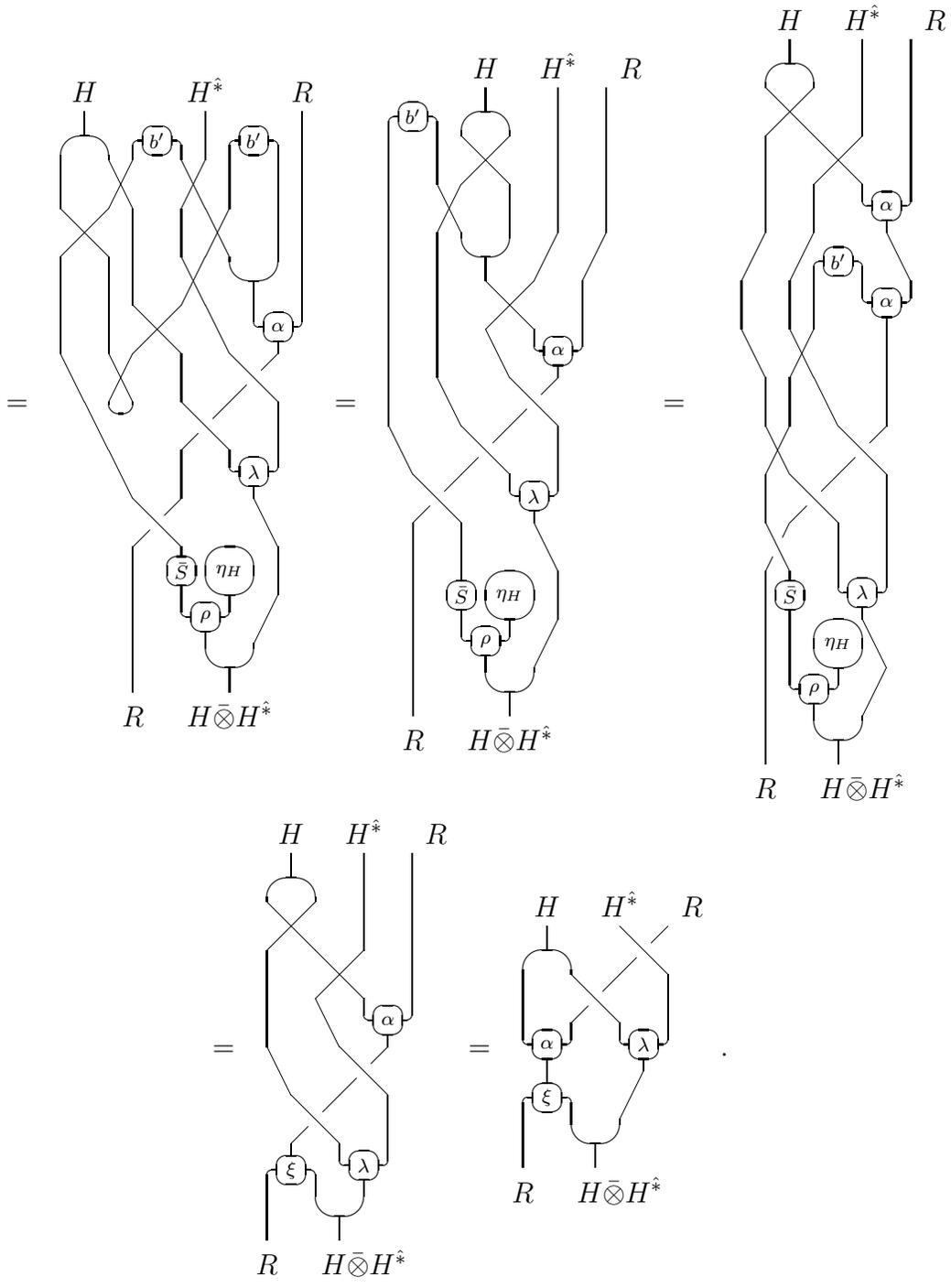
\[
\begin{tangle}
\object{H}\step[2]\object{H^{\hat{*}}}\Step\object{R}\\
\tu \lambda\step\td \xi\\
\step\x\Step\id\\
\step\id\Step\cu\\
\object{R}\step[5]\object{H\bar{\otimes} H^{\hat{*}}}\\
\end{tangle}
\step=\step
\begin{tangle}
\step\object{H}\step[2]\object{H^{\hat{*}}}\step[3]\object{R}\\
\cd\Step\x\\
\id\Step\x\Step\id\\
\cu\Step\tu \lambda\\
\td \xi\Step\dd\\
\id\Step\cu\\
\object{R}\step[4]\object{H\bar{\otimes} H^{\hat{*}}}\\
\end{tangle}\ \ \ .
\step[5] ......(*)
\]

\[
\hbox{the left side}\step=\step
\begin{tangle}
\object{H}\step[3]\object{H^{\hat{*}}}\step[4]\object{R}\step[3]\object{\eta_H}\\
\tu \lambda\step\ro {b'}\Step\id\step[3]\id\\
\step\id\Step\O {\bar{S}}\Step\tu \alpha\step[3]\id\\
\step\id\Step\d\Step\id\step[3]\dd\\
\step\d\Step\x\Step\dd\\
\Step\x\Step\tu \rho\\
\Step\id\Step\d\Step\id\\
\Step\id\step[3]\cu\\
\Step\object{R}\step[5]\object{H\bar{\otimes} H^{\hat{*}}}\\
\end{tangle}
\step\stackrel{\hbox{by Lemma1.5(1)}}{\step[3]=\step[3]}\step
\begin{tangle}
\object{H}\step[3]\object{H^{\hat{*}}}\step[5]\object{R}\step[3]\object{\eta_H}\\
\id\Step\id\Step\ro {b'}\Step\id\step[3]\id\\
\id\Step\id\Step\O {\bar{S}}\Step\tu \alpha\step[3]\id\\
\id\Step\id\Step\id\Step\dd\step[4]\id\\
\id\Step\id\Step\x\step[4]\dd\\
\id\Step\x\Step\id\step[3]\dd\\
\x\Step\XX\step[2]\dd\\
\id\Step\XX\Step\XX\\
\id\step\cd\step\XX\step\cd\\
\id\step\O S\Step\id\step\id\Step\id\step\XX\\
\id\step\XX\step\id\Step\X\Step\id\\
\id\step\id\Step\X\Step\id\step\id\Step\id\\
\id\step\id\step\dd\step\XX\step\d\step\id\\
\id\step\id\step\XX\Step\XX\step\id\\
\id\step\id\step\tu \rightharpoonup\Step\tu
\leftharpoonup\step\id\\
\id\step\tu \rho\step[4]\tu \lambda\\
\id\Step\d\step[4]\dd\\
\id\step[3]\Cu\\
\object{R}\step[5]\object{H\bar{\otimes} H^{\hat{*}}}\\
\end{tangle} \]
\[
\step=\step
\begin{tangle}
\object{H}\step[3]\object{H^{\hat{*}}}\step[5]\object{R}\step[3]\object{\eta_H}\\
\id\Step\id\Step\ro {b'}\Step\id\step[3]\id\\
\id\Step\id\Step\O {\bar{S}}\Step\tu \alpha\step[3]\id\\
\id\Step\id\Step\id\Step\dd\step[4]\id\\
\id\Step\id\Step\x\step[4]\dd\\
\id\Step\x\Step\id\step[3]\dd\\
\x\Step\XX\step[2]\dd\\
\id\Step\XX\Step\XX\\
\id\step\cd\step\XX\Step\id\\
\id\step\O {\bar S}\Step\id\step\id\Step\id\Step\id\\
\id\step\XX\step\id\Step\id\Step\id\\
\id\step\id\Step\X\Step\id\Step\id\\
\id\step\tu \rho\step\XX\Step\id\\
\id\Step\id\Step\tu \leftharpoonup\step\dd\\
\id\Step\id\step[3]\tu \lambda\\
\id\Step\Cu\\
\object{R}\step[5]\object{H\bar{\otimes} H^{\hat{*}}}\\
\end{tangle}
\step=\step
\begin{tangle}
\object{H}\step[3]\object{H^{\hat{*}}}\step[5]\object{R}\\
\id\Step\id\Step\ro {b'}\Step\id\\
\id\Step\XX\Step\tu \alpha\\
\XX\Step\id\Step\dd\\
\id\Step\id\Step\x\\
\id\Step\x\Step\id\\
\x\Step\id\Step\id\\
\id\step\cd\step\d\step\id\\
\id\step\O {\bar S}\Step\XX\step\id\\
\id\step\id\Step\tu \leftharpoonup\step\id\\
\id\step\id\step[3]\tu \lambda\\
\id\step\id\step\obox 2{\eta_H}\step\id\\
\id\step\tu \rho\step\dd\\
\id\Step\cu\\
\object{R}\step[4]\object{H\bar{\otimes} H^{\hat{*}}}\\
\end{tangle}
\]
\[
\step=\step
\begin{tangle}
\step\object{H}\step[5]\object{H^{\hat{*}}}\step[4]\object{R}\\
\cd\step\ro {b'}\step\id\step\ro {b'}\step\id\\
\id\Step\X\Step\X\step\id\Step\id\step\id\\
\XX\step\id\Step\id\step\X\Step\id\step\id\\
\id\Step\id\step\id\Step\X\step\cu\step\id\\
\id\Step\id\step\XX\step\d\step\tu \alpha\\
\d\step\X\Step\id\Step\x\\
\step\d\hev\Step\x\Step\id\\
\Step\d\Step\id\Step\tu \lambda\\
\step[3]\x\step[3]\d\\
\step[3]\id\Step\O {\bar{S}}\step\obox 2{\eta_H}\step\id\\
\step[3]\id\Step\tu \rho\step\dd\\
\step[3]\id\step[3]\cu\\
\step[3]\object{R}\step[4]\object{H\bar{\otimes} H^{\hat{*}}}\\
\end{tangle}
\step=\step
\begin{tangle}
\step[4]\object{H}\step[3]\object{H^{\hat{*}}}\step[3]\object{R}\\
\ro {b'}\step\cd\step[2]\id\Step\id\\
\id\Step\id\step\XX\Step\id\Step\id\\
\id\Step\X\Step\id\Step\id\Step\id\\
\id\Step\id\step\cu\step\dd\step\dd\\
\id\Step\id\step [2]\XX\step [2]\id\\
\id\Step\id\step [2]\d\step\tu \alpha\\
\id\Step\d\Step\x\\
\d\Step\x\Step\id\\
\step\x\Step\tu \lambda\\
\step\id\Step\id\step[3]\d\\
\step\id\Step\O {\bar{S}}\step\obox 2{\eta_H}\step\id\\
\step\id\Step\tu \rho\step\dd\\
\step\id\step[3]\cu\\
\step\object{R}\step[4]\object{H\bar{\otimes} H^{\hat{*}}}\\
\end{tangle}
\step=\step
\begin{tangle}
\step[3]\object{H}\step[3]\object{H^{\hat{*}}}\step[3]\object{R}\\
\Step\cd\Step\id\Step\id\\
\Step\XX\Step\id\Step\id\\
\step [2]\id \step[2]\XX  \Step\id\\
\step [2]\id \step[2]\id \step [2]\tu \alpha\\
\step \dd \step \dd\ro {b'}\step [1]\d\\
\step \id \step[2]\id \step \id\step [2]\tu \alpha\\
\step \d \step \X\step[3]\id\\
\step[2]\id \step[1]\id \step \d \step [2] \id \\
\step[2]\X\step[2]\x \\
\step [2] \id \step [1] \x \step [2]\id \\
\step[2]\hx\step[2]\id\Step\id\\
\step [2]\id \step\O {\bar S}\step \step \tu \lambda\\
\step [2]\id\step\id\step\obox 2{\eta_H}\d\\
\step [2]\id\step\tu \rho\step\dd\\
\step [2]\id\step[2]\cu\\
\step [2]\object{R}\step[4]\object{H\bar{\otimes} H^{\hat{*}}}\\
\end{tangle}\]
\[
\step=\step
\begin{tangle}
\step\object{H}\step[3]\object{H^{\hat{*}}}\step[3]\object{R}\\
\cd\Step\id\Step\id\\
\XX\Step\id\Step\id\\
\id\Step\XX\Step\id\\
\id\Step\d\step\tu \alpha\\
\d\Step\x\\
\step\x\Step\id\\
\td \xi\step\tu \lambda\\
\id\Step\cu\\
\object{R}\step[4]\object{H\bar{\otimes} H^{\hat{*}}}\\
\end{tangle}
\step=\step
\begin{tangle}
\step\object{H}\step[3]\object{H^{\hat{*}}}\step[3]\object{R}\\
\cd\Step\x\\
\id\Step\x\Step\id\\
\tu \alpha\Step\tu \lambda\\
\td \xi\Step\dd\\
\id\Step\cu\\
\object{R}\step[4]\object{H\bar{\otimes} H^{\hat{*}}}\\
\end{tangle}\ \ \ .
\]
Thus relation (*) holds.

Next, we check that $\xi$ is algebraic. We see that
\[
\begin{tangle}
\object{R}\Step\object{R}\\
\O \xi\Step\O \xi\\
\cu\\
\step\object{R\otimes(H\bar{\otimes}H^{\hat{*}})}\\
\end{tangle}
\step[3]=\step
\begin{tangle}
\step\object{R}\step[6]\object{R}\\
\td \phi\step[4]\td \phi\\
\O {\bar{S}}\Step\id\step[4]\O {\bar{S}}\Step\id\\
\x\step\obox 2{\eta_H}\step\x\step\obox 2{\eta_H}\\
\id\Step\tu \rho\step\dd\Step\tu \rho\\
\d\Step\x\step[4]\id\\
\step\cu\Step\Cu\\
\Step\object{R}\step[5]\object{H\bar{\otimes} H^{\hat{*}}}\\
\end{tangle}
\step=\step
\begin{tangle}
\step\object{R}\step[3]\object{R}\\
\td \phi\step\td \phi\\
\O {\bar{S}}\Step\id\step\O {\bar{S}}\Step\id\\
\x\step\id\Step\id\\
\id\Step\X\Step\id\\
\id\Step\id\step\x\\
\id\Step\hx\Step\id\\
\cu\step\cu\\\
\step\id\step[3]\d\step\obox 2{\eta_H}\\
\step\id\step[4]\tu \rho\\
\step\object{R}\step[5]\object{H\bar{\otimes} H^{\hat{*}}}\\
\end{tangle}
\step=\step[4]
\begin{tangle}
\object{R}\Step\object{R}\\
\cu\\
\step\O \xi\\
\step\object{R\otimes(H\bar{\otimes}H^{\hat{*}})}\\
\end{tangle}
\]
and obviously
\[
\begin{tangle}
\object{\eta_R}\\
\id\\
\O \xi\\
\id\\
\object{R\otimes(H\bar{\otimes}H^{\hat{*}})}\\
\end{tangle}
\step[6]=\step[6]
\begin{tangle}
\object{\eta_{R\otimes(H\bar{\otimes}H^{\hat{*}})}}\\
\id\\
\id\\
\id\\
\object{R\otimes(H\bar{\otimes}H^{\hat{*}})}\\
\end{tangle}\ \ \ \ \ \ \ \ \ .
\]
Thus $\xi$ is algebraic.

Now we show that $\Phi'$ is algebraic.
\[
\begin{tangle}
\object{(R\# H)\# H^{\hat{*}}}\step[11]\object{(R\# H)\# H^{\hat{*}}}\\
\step\nw2\step[6]\ne2\\
\step[2]\obox 2{\Phi'}\step[2]\obox 2{\Phi'}\\
\step[3]\Cu\\
\step[6]\object{R\otimes(H\bar{\otimes}H^{\hat{*}})}\\
\end{tangle}
\step[3]=\Step
\begin{tangle}
\step\object{R}\step[4]\object{H}\step[3]\object{H^{\hat{*}}}
\step[3]\object{R}\step[4]\object{H}\step[3]\object{H^{\hat{*}}}\\
\td \xi\step[3]\id\Step\id\step[3]\td \xi\step[3]\id\Step\id\\
\id\Step\id\step[3]\tu \lambda\Step\dd\Step\id\step[3]\tu \lambda\\
\id\Step\Cu\Step\dd\step[3]\Cu\\
\id\step[4]\d\Step\dd\step[5]\ne2\\
\d\step[4]\x\step[4]\dd\\
\step\Cu\Step\Cu\\
\step[3]\object{R}\step[6]\object{H\bar{\otimes}H^{\hat{*}}}\\
\end{tangle}\]
\[
\step=\step
\begin{tangle}
\step\object{R}\step[3]\object{H}\step[2]\object{H^{\hat{*}}}
\step[3]\object{R}\step[3]\object{H}\step[2]\object{H^{\hat{*}}}\\
\td \xi\step[2]\id\Step\id\step[2]\td \xi\step[2]\id\Step\id\\
\id\Step\id\Step\tu \lambda\step\dd\step\dd\Step\tu \lambda\\
\id\Step\id\step[3]\x\Step\id\step[4]\id\\
\id\Step\id\Step\dd\Step\cu\step[3]\dd\\
\id\Step\x\step[4]\Cu\\
\cu\Step\d\step[4]\dd\\
\step\id\step[4]\Cu\\
\step\object{R}\step[6]\object{H\bar{\otimes}H^{\hat{*}}}\\
\end{tangle}
\stackrel{\hbox{by (*)}}{=}
\begin{tangle}
\step\object{R}\step[3]\object{H}\step[2]\object{H^{\hat{*}}}
\step[2]\object{R}\step[2]\object{H}\step[3]\object{H^{\hat{*}}}\\
\td \xi\step\cd\step\x\Step\id\Step\id\\
\id\Step\id\step\id\Step\hx\Step\id\Step\tu \lambda\\
\id\Step\id\step\tu \alpha\step\tu \lambda\step[3]\id\\
\id\Step\id\step\td \xi\Step\id\step[4]\id\\
\id\Step\hx\Step\cu\step[3]\dd\\
\cu\step\d\Step\Cu\\
\step\id\step[3]\Cu\\
\step\object{R}\step[5]\object{H\bar{\otimes}H^{\hat{*}}}\\
\end{tangle}
\]
\[
\step\stackrel{\hbox{by Lemma1.4}}{=}\step
\begin{tangle}
\step\object{R}\step[3]\object{H}\step[2]\object{H^{\hat{*}}}
\step[2]\object{R}\step[2]\object{H}\step[3]\object{H^{\hat{*}}}\\
\td \xi\step\cd\step\x\Step\id\Step\id\\
\id\Step\id\step\id\Step\hx\step\cd\step\id\Step\id\\
\id\Step\id\step\tu \alpha\step\id\step\id\Step\X\Step\id\\
\id\Step\id\step\td \xi\step\id\step\tu \rightharpoonup\step\cu\\
\id\Step\hx\Step\id\step\cu\step[2]\ne2\\
\cu\step\id\Step\d\step\tu \lambda\\
\step\id\Step\id\step[3]\cu\\
\step\id\Step\Cu\\
\step\object{R}\step[5]\object{H\bar{\otimes}H^{\hat{*}}}\\
\end{tangle}
\step\stackrel{\hbox{ since } \xi \hbox{ is algebraic }}{=}\step
\begin{tangle}
\object{R}\step[2]\object{H}\step[2]\object{H^{\hat{*}}}
\step[2]\object{R}\step[2]\object{H}\step[3]\object{H^{\hat{*}}}\\
\id\step\cd\step\x\Step\id\Step\id\\
\id\step\id\Step\hx\step\cd\step\id\Step\id\\
\id\step\tu \alpha\step\id\step\id\Step\X\Step\id\\
\cu\Step\id\step\tu \rightharpoonup\step\cu\\
\td \xi\Step\cu\Step\ne2\\
\id\Step\id\step[3]\tu \lambda\\
\id\Step\Cu\\
\object{R}\step[5]\object{H\bar{\otimes}H^{\hat{*}}}\\
\end{tangle} \]
\[
\step=\step[6]
\begin{tangle}
\object{(R\# H)\# H^{\hat{*}}}\step[11]\object{(R\# H)\# H^{\hat{*}}}\\
\step\nw2\step[6]\ne2\\
\step[3]\Cu\\
\step[4]\obox 2{\Phi'}\\
\step[5]\id\\
\step[4]\object{R\otimes(H\bar{\otimes}H^{\hat{*}})}\\
\end{tangle}\ \ \ \ \ \ \ \ \ \ .
\]

It is clear that
\[
\begin{tangle}
\object{\eta_{(R\# H)\# H^{\hat{*}}}}\\
\step\id\\
\obox 2{\Phi'}\\
\step\id\\
\step\object{R\otimes(H\bar{\otimes}H^{\hat{*}})}\\
\end{tangle}
\step[6]=\step[6]
\begin{tangle}
\object{\eta_{R\otimes(H\bar{\otimes}H^{\hat{*}})}}\\
\id\\
\id\\
\id\\
\object{R\otimes(H\bar{\otimes}H^{\hat{*}})}\\
\end{tangle}\ \ \ \ \ \ \ \ \ \ \ \ \ .\] Thus $\Phi '$ is algebraic. \ \
\begin{picture}(5,5)
\put(0,0){\line(0,1){5}}\put(5,5){\line(0,-1){5}}
\put(0,0){\line(1,0){5}}\put(5,5){\line(-1,0){5}}
\end{picture}\\

The proof of Lemma3.2.\\
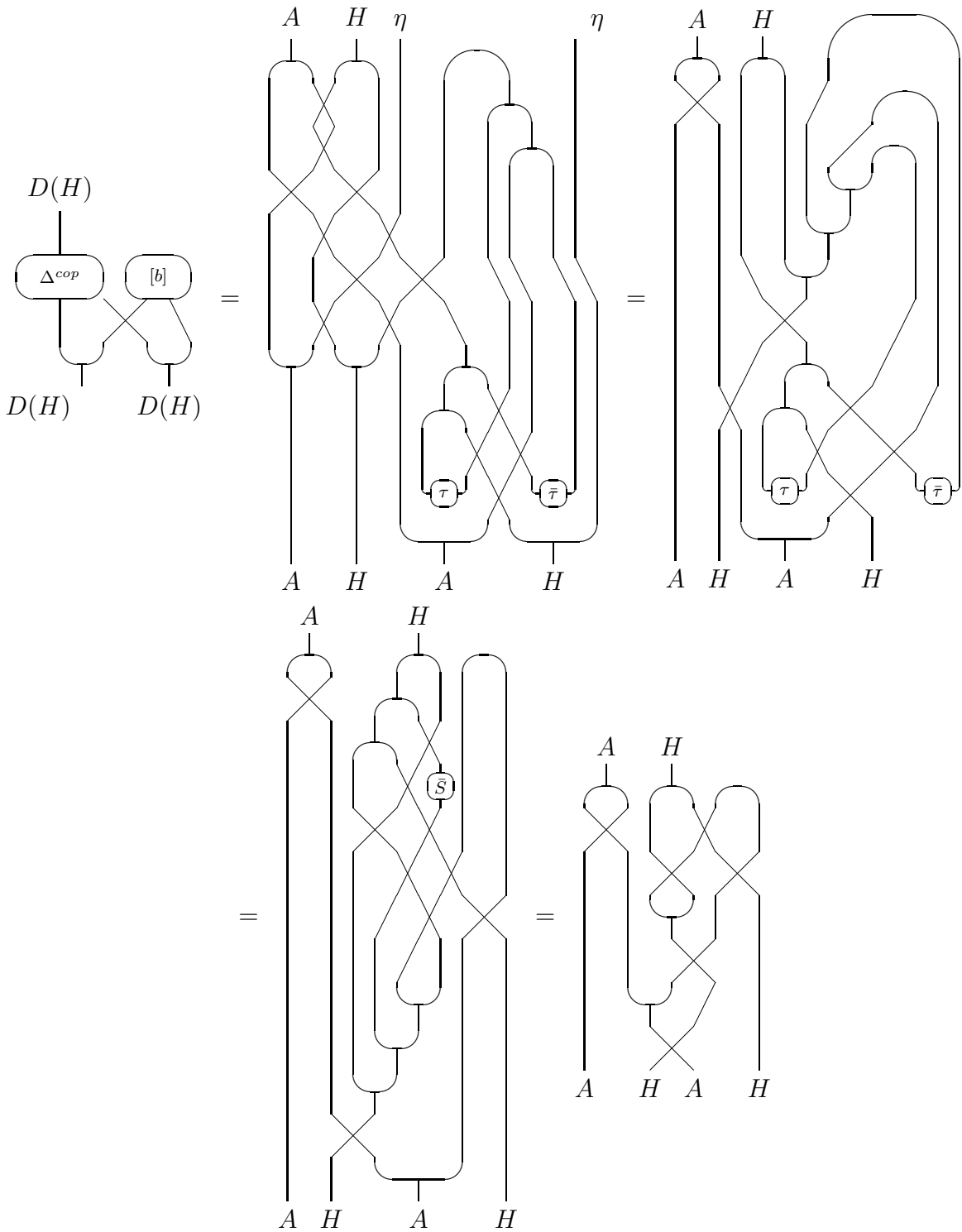
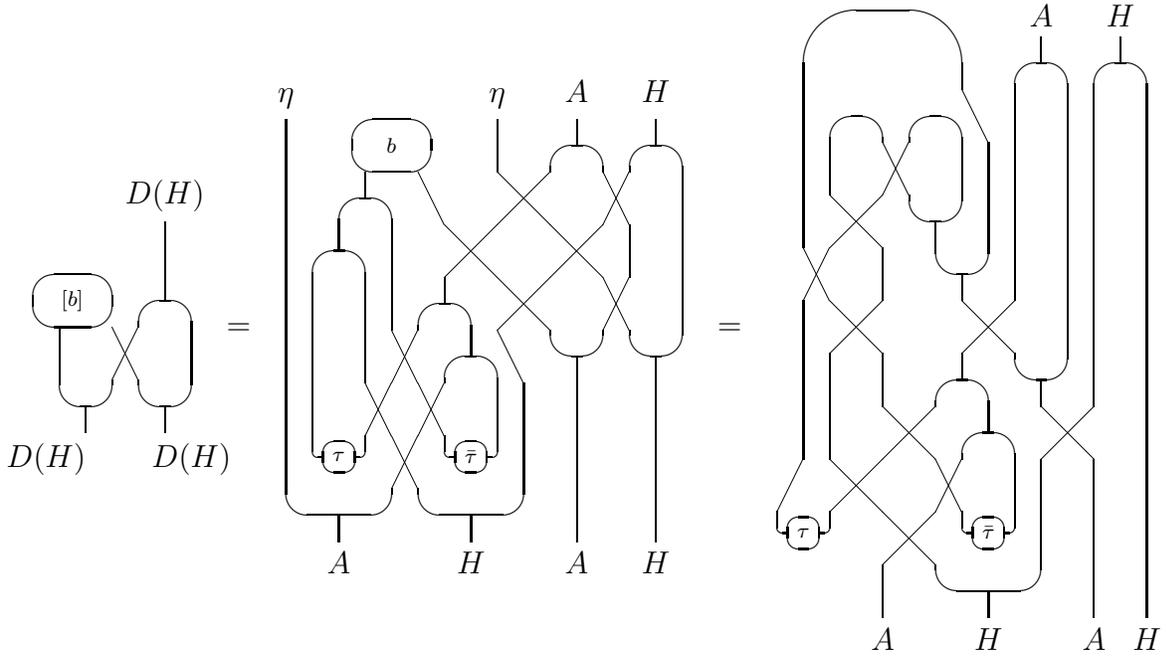
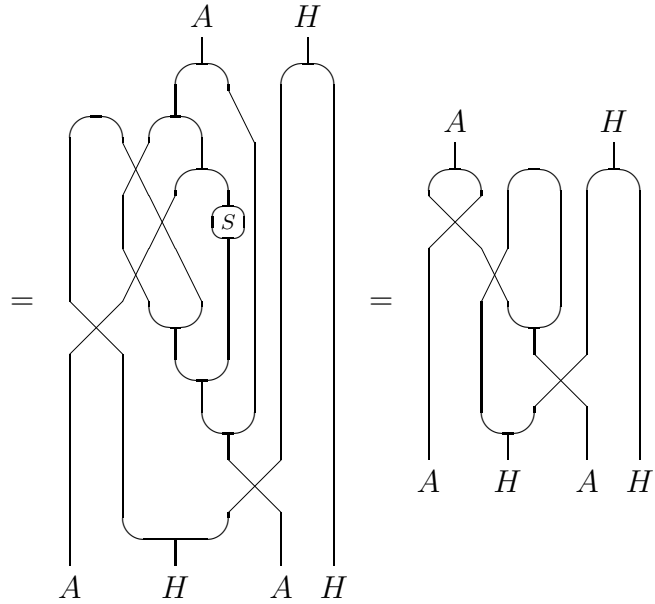
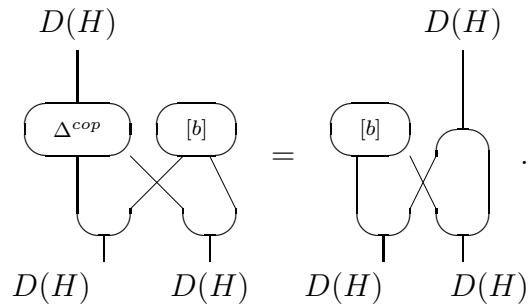
\[
\begin{tangle}
\step[2]\object{D(H)}\\
\step[2]\id\\
\obox 4{\Delta^{cop}}\step\obox 3{[b]}\\
\Step\id\step[2]\XX\step\d\\
\Step\cu\Step\cu\\
\step\object{D(H)}\step[6]\object{D(H)}\\
\end{tangle}
\step=\step
\begin{tangle}
\step\object{A}\step[3]\object{H}\Step\object{\eta}\step[9]\object{\eta}\\
\cd\step\cd\step\id\Step\Coev\step[5]\id\\
\id\Step\X\Step\id\step\id\Step\id\step[2]\cd\Step\id\\
\id\Step\X\Step\id\step\id\Step\id\Step\id\step\cd\step\id\\
\XX\step\XX\step\id\Step\id\Step\id\step\id\Step\id\step\id\\
\id\Step\X\Step\X\Step\id\Step\id\step\id\Step\id\step\id\\
\id\Step\id\step\XX\step\XX\Step\d\d\step\d\d\\
\id\Step\X\Step\X\Step\d\Step\id\step\id\Step\id\step\id\\
\cu\step\cu\step\id\step[2]\cd\step\id\step\id\Step\id\step\id\\
\step\id\step[3]\id\Step\id\step\cd\step\X\step\id\Step\id\step\id\\
\step\id\step[3]\id\Step\id\step\id\Step\X\step\X\Step\id\step\id\\
\step\id\step[3]\id\Step\id\step\coro \tau \step\X\step\coro {\bar \tau }\step\id\\
\step\id\step[3]\id\Step\Cu\step\Cu\\
\step\object{A}\step[3]\object{H}\step[4]\object{A}\step[5]\object{H}\\
\end{tangle}
\step=\step
\begin{tangle}
\step\object{A}\step[3]\object{H}\\
\cd\step\cd\Step\COEV\\
\XX\step\id\Step\id\step\dd\Step\Coev\step[3]\id\\
\id\Step\id\step\id\Step\id\step\id\Step\ne2\coev\step\id\step\id\\
\id\Step\id\step\id\Step\id\step\id\step\cu\Step\id\step\id\step\id\\
\id\Step\id\step\id\Step\id\step\cu\step[3]\id\step\id\step\id\\
\id\Step\id\step\d\step\cu\step[4]\id\step\id\step\id\\
\id\Step\id\Step\XX\step[4]\dd\step\id\step\id\\
\id\Step\id\step\dd\step\cd\Step\dd\Step\id\step\id\\
\id\Step\X\step\cd\step\XX\Step\dd\step\id\\
\id\Step\id\step\id\step\id\Step\X\Step\XX\Step\id\\
\id\Step\id\step\id\step\coro \tau \step\XX\Step\coro {\bar \tau }\\
\id\Step\id\step\Cu\Step\id\\
\object{A}\Step\object{H}\step[3]\object{A}\step[4]\object{H}\\
\end{tangle} \]
\[
\step=\step
\begin{tangle}
\step\object{A}\step[5]\object{H}\\
\cd\step[3]\cd\step\coev\\
\XX\Step\cd\step\id\step\id\Step\id\\
\id\Step\id\step\cd\step\X\step\id\Step\id\\
\id\Step\id\step\id\Step\X\step\O {\bar{S}}\step\id\Step\id\\
\id\Step\id\step\XX\step\X\step\id\Step\id\\
\id\Step\id\step\id\Step\X\step\X\Step\id\\
\id\Step\id\step\id\step\dd\step\X\step\XX\\
\id\Step\id\step\id\step\id\step\dd\step\id\step\id\Step\id\\
\id\Step\id\step\id\step\id\step\cu\step\id\Step\id\\
\id\Step\id\step\id\step\cu\Step\id\Step\id\\
\id\Step\id\step\cu\step[3]\id\Step\id\\
\id\Step\XX\step[4]\id\Step\id\\
\id\Step\id\Step\Cu\Step\id\\
\object{A}\Step\object{H}\step[4]\object{A}\step[4]\object{H}\\
\end{tangle}
\step=\step
\begin{tangle}
\step\object{A}\step[3]\object{H}\\
\cd\step\cd\step\coev\\
\XX\step\id\Step\X\Step\id\\
\id\Step\id\step\XX\step\XX\\
\id\Step\id\step\cu\step\id\Step\id\\
\id\Step\id\Step\XX\Step\id\\
\id\Step\cu\step\dd\Step\id\\
\id\step[3]\XX\step[3]\id\\
\object{A}\step[3]\object{H}\Step\object{A}\step[3]\object{H}\\
\end{tangle}
\]

and
\[
\begin{tangle}
\step[5]\object{D(H)}\\
\step[5]\id\\
\obox 3{[b]}\step\cd\\
\step\id\step[2]\X\Step\id\\
\step\cu\step\cu\\
\step[0.5]\object{D(H)}\step[5.5]\object{D(H)}\\
\end{tangle}
\step=\step
\begin{tangle}
\object{\eta}\step[8]\object{\eta}\step[3]\object{A}\step[3]\object{H}\\
\id\step[2.5]\obox 3{b}\step[2.5]\id\Step\cd\step\cd\\
\id\Step\cd\step\d\Step\XX\Step\X\Step\id\\
\id\step\cd\step\id\Step\XX\Step\XX\step\id\Step\id\\
\id\step\id\Step\id\step\id\step\cd\step\XX\Step\X\Step\id\\
\id\step\id\Step\id\step\X\step\cd\d\step\cu\step\cu\\
\id\step\id\Step\X\step\X\Step\id\step\id\Step\id\step[3]\id\\
\id\step\coro \tau \step\X\step\coro {\bar \tau }\step\id\Step\id\step[3]\id\\
\Cu\step\Cu\Step\id\step[3]\id\\
\Step\object{A}\step[5]\object{H}\step[4]\object{A}\step[3]\object{H}\\
\end{tangle}
\step=\step
\begin{tangle}
\step \step[9]\object{A}\step[3]\object{H}\\
\step \COEV\step[6]\cd\step\cd\\
\step \id\step\coev\step\coev\d\step\id\Step\id\step\id\Step\id\\
\step \id\step\id\Step\X\Step\id\step\id\step\id\Step\id\step\id\Step\id\\
\step \id\step\XX\step\cu\step\id\step\id\Step\id\step\id\Step\id\\
\step \X\Step\id\Step\cu\step\id\Step\id\step\id\Step\id\\
\step \id\step\XX\step[3]\XX\Step\id\step\id\Step\id\\
\step \id\step\id\Step\id\Step\cd\step\cu\step\id\Step\id\\
\step \id\step\id\Step\XX\step\cd\step\XX\Step\id\\
\dd\step\XX\Step\X\Step\id\step\id\Step\id\Step\id\\
\coro \tau \Step\XX\step\coro {\bar \tau }\step\id\Step\id\Step\id\\
\step \step[3]\id\Step\Cu\Step\id\Step\id\\
\step \step[3]\object{A}\step[4]\object{H}\step[4]\object{A}\Step\object{H}\\
\end{tangle}\]
\[
\step=\step
\begin{tangle}
\step[5]\object{A}\step[4]\object{H}\\
\step[4]\cd\Step\cd\\
\coev\step\cd\step\d\step\id\Step\id\\
\id\Step\X\step\cd\step\id\step\id\Step\id\\
\id\Step\id\step\X\Step\O S\step\id\step\id\Step\id\\
\id\Step\X\step\d\step\id\step\id\step\id\Step\id\\
\XX\step\cu\step\id\step\id\step\id\Step\id\\
\id\Step\id\Step\cu\step\id\step\id\Step\id\\
\id\Step\id\step[3]\cu\step\id\Step\id\\
\id\Step\id\step[4]\XX\Step\id\\
\id\Step\Cu\Step\id\Step\id\\
\object{A}\step[4]\object{H}\step[4]\object{A}\Step\object{H}\\
\end{tangle}
\step=\step
\begin{tangle}
\step\object{A}\step[6]\object{H}\\
\cd\step\coev\step\cd\\
\XX\step\id\Step\id\step\id\Step\id\\
\id\Step\X\Step\id\step\id\Step\id\\
\id\Step\id\step\cu\step\id\Step\id\\
\id\Step\id\Step\XX\Step\id\\
\id\Step\cu\Step\id\Step\id\\
\object{A}\step[3]\object{H}\step[3]\object{A}\step[2]\object{H}\\
\end{tangle}\ \ \ .
\]

Thus
\[
\begin{tangle}
\step[2]\object{D(H)}\\
\step[2]\id\\
\obox 4{\Delta^{cop}}\step\obox 3{[b]}\\
\Step\id\step[2]\XX\step\d\\
\Step\cu\Step\cu\\
\step\object{D(H)}\step[6]\object{D(H)}\\
\end{tangle}
\step=\step
\begin{tangle}
\step[5]\object{D(H)}\\
\step[5]\id\\
\obox 3{[b]}\step\cd\\
\step\id\step[2]\X\Step\id\\
\step\cu\step\cu\\
\step[0.5]\object{D(H)}\step[5.5]\object{D(H)}\\
\end{tangle}\ \ \ .
\]
We complete the proof.\ \
\begin{picture}(5,5)
\put(0,0){\line(0,1){5}}\put(5,5){\line(0,-1){5}}
\put(0,0){\line(1,0){5}}\put(5,5){\line(-1,0){5}}
\end{picture}\\

{\bf Remark:} If $U$ and $V$ have left dual $U^*$ and $V^*$,
respectively, then $U^*\otimes V^*$ and $V^* \otimes U^*$ both are
the left duals of $U\otimes V$. Their evaluations and
coevaluations are $$ d_{U\otimes V}=(d_{U} \otimes d_V) (id _{U^*}
\otimes C_{V^*, U}\otimes  id _V), \ \ \ \ b_{U\otimes V}=(id
_U\otimes (C_{V, U^*})^{-1}\otimes id_ {V^*})(b_U \otimes b_V);$$\
$$ d_{U\otimes V}=d_V (id _{V^*}\otimes d_U \otimes id _V ) ,
\ \ \ \ b_{U\otimes V}=(id _{U^*}\otimes b_V \otimes id_U)b_U,
$$ respectively. In this paper, we use the second.  $H^{\hat * }$ is the left dual of $H$ under the first.

\end{document}